\documentclass[12pt]{article}
\usepackage[]{amsmath,amssymb}
\usepackage{latexsym}
\usepackage{cite}

\newtheorem{definition}{Definition}[section]
\newtheorem{theorem}[definition]{Theorem}
\newtheorem{lemma}[definition]{Lemma}

\newtheorem{conjecture}[definition]{Conjecture}

\newtheorem{note}[definition]{Note}
\newtheorem{assumption}[definition]{Assumption}
\newtheorem{proposition}[definition]{Proposition}

\def\F{\mathbb F}
\def\K{\mathbb F}

\begin{document}
\title{\bf Tridiagonal pairs of $q$-Racah type}
\author{
Tatsuro Ito\footnote{Supported in part by JSPS grant
18340022.} $\;$   and
Paul Terwilliger\footnote{This author gratefully acknowledges 
support from the FY2007 JSPS Invitation Fellowship Program
for Reseach in Japan (Long-Term), grant L-07512.}
}
\date{}

\maketitle
\begin{abstract}
Let $\K$ denote an algebraically closed
 field and let $V$ denote a vector space over $\K$ with
finite positive dimension.
We consider a pair of linear transformations $A:V \to V$
and $A^*:V \to V$ that satisfy the following conditions:
(i)
each of $A,A^*$ is diagonalizable;
(ii)
there exists an ordering $\lbrace V_i\rbrace_{i=0}^d$ of the eigenspaces of
$A$ such that
$A^* V_i \subseteq V_{i-1} + V_{i} + V_{i+1}$ for $0 \leq i \leq d$,
where $V_{-1}=0$ and $V_{d+1}=0$;
(iii)
there exists an ordering $\lbrace V^*_i\rbrace_{i=0}^\delta$
of the eigenspaces of $A^*$ such that
$A V^*_i \subseteq V^*_{i-1} + V^*_{i} + V^*_{i+1}$ for
 $0 \leq i \leq \delta$,
where $V^*_{-1}=0$ and $V^*_{\delta+1}=0$;
(iv)
there is no subspace $W$ of $V$ such that
$AW \subseteq W$, $A^* W \subseteq W$, $W \neq 0$, $W \neq V$.
We call such a pair a {\it tridiagonal pair} on $V$.
It is known that $d=\delta$.
For $0 \leq i \leq d$ let $\theta_i$ (resp. $\theta^*_i$)
denote the eigenvalue of $A$ (resp. $A^*$) associated with
$V_i$ (resp. $V^*_i$). The pair $A,A^*$ is said to have
{\it $q$-Racah type} whenever 
$\theta_i = a + bq^{2i-d}+ c q^{d-2i}$ and
$\theta^*_i = a^* + b^*q^{2i-d}+c^*q^{d-2i}$ for $0 \leq i \leq d$,
where $q, 
a,b,c,a^*,b^*,c^*$
 are scalars in $\K$
with $q,b,c,b^*,c^*$ nonzero and $q^2 \not\in \lbrace
1,-1\rbrace$.
This type is the most general one.
We classify up to isomorphism the tridiagonal pairs
over $\K$ that have $q$-Racah type.  Our proof
involves the representation theory of
the quantum affine algebra
 $U_q(\widehat{ \mathfrak{sl}}_2)$.

 \end{abstract}

\section{Tridiagonal pairs}

\noindent 
Throughout this paper $\K$ denotes a field
and $\overline \K$ denotes the algebraic closure of $\K$.

\medskip
\noindent 
We begin by recalling the notion of a tridiagonal pair. 
We will use the following terms.
Let $V$ denote a vector space over $\K$ with finite
positive dimension.
For a 
 linear transformation $A:V\to V$
and a
subspace $W \subseteq V$,
we call $W$ an
 {\it eigenspace} of $A$ whenever 
 $W\not=0$ and there exists $\theta \in \K$ such that 
$W=\lbrace v \in V \;\vert \;Av = \theta v\rbrace$;
in this case $\theta$ is the {\it eigenvalue} of
$A$ associated with $W$.
We say that $A$ is {\it diagonalizable} whenever
$V$ is spanned by the eigenspaces of $A$.

\begin{definition}  
{\rm \cite[Definition~1.1]{TD00}}
\label{def:tdp}
\rm
Let $V$ denote a vector space over $\K$ with finite
positive dimension. 
By a {\it tridiagonal pair} (or {\it $TD$ pair})
on $V$
we mean an ordered pair of linear transformations
$A:V \to V$ and 
$A^*:V \to V$ 
that satisfy the following four conditions.
\begin{enumerate}
\item Each of $A,A^*$ is diagonalizable.
\item There exists an ordering $\lbrace V_i\rbrace_{i=0}^d$ of the  
eigenspaces of $A$ such that 
\begin{equation}
A^* V_i \subseteq V_{i-1} + V_i+ V_{i+1} \qquad \qquad (0 \leq i \leq d),
\label{eq:t1}
\end{equation}
where $V_{-1} = 0$ and $V_{d+1}= 0$.
\item There exists an ordering $\lbrace V^*_i\rbrace_{i=0}^{\delta}$ of
the  
eigenspaces of $A^*$ such that 
\begin{equation}
A V^*_i \subseteq V^*_{i-1} + V^*_i+ V^*_{i+1} 
\qquad \qquad (0 \leq i \leq \delta),
\label{eq:t2}
\end{equation}
where $V^*_{-1} = 0$ and $V^*_{\delta+1}= 0$.
\item There does not exist a subspace $W$ of $V$ such  that $AW\subseteq W$,
$A^*W\subseteq W$, $W\not=0$, $W\not=V$.
\end{enumerate}
We say the pair $A,A^*$ is {\it over $\K$}.
We call $V$ the 
{\it underlying
 vector space}.
\end{definition}

\begin{note}
\label{lem:convention}
\rm
According to a common notational convention $A^*$ denotes 
the conjugate-transpose of $A$. We are not using this convention.
In a TD pair $A,A^*$ the linear transformations $A$ and $A^*$
are arbitrary subject to (i)--(iv) above.
\end{note}

\medskip
\noindent 
We now give some background on TD pairs;
for more information we refer the reader to
\cite{
TD00,
shape,
tdanduq,
NN,
qtet,
Ev,
madrid}.
The concept of a TD pair
originated in the study of the 
($P$ and $Q$)-polynomial association schemes
\cite{BanIto}
and their relationship to the Askey scheme
of orthogonal polynomials
\cite{AWil,
KoeSwa}.
The concept is implicit in
\cite[p.~263]{BanIto},
\cite{Leon}
and more explicit in 
\cite[Theorem~2.1]{TersubI}.
A systematic study began in
\cite{TD00}.
As research progressed, connections were found to
representation theory
\cite{
bas1,hasan2,
neubauer,
Ha,
tdanduq,
qtet,
IT:Krawt,
IT:aug,
Koelink3,
cite37,
cite40,
noumi1,  
Rosengren,
Rosengren2,
  qSerre,
aw},
partially ordered sets
    \cite{lsint},
the bispectral problem 
\cite{GYLZmut, GH7,GH1,Zhidd},
statistical mechanical models
 \cite{bas1,bas2,bas3,bas4,bas5,bas6,bas7,DateRoan2, Davfirst, Da,Onsager},
and classical mechanics
\cite{LPcm}.

\medskip
\noindent 
We now recall some basic facts about TD pairs.
Let $A,A^*$ denote a TD pair
on $V$, as in Definition 
\ref{def:tdp}. By
\cite[Lemma 4.5]{TD00}
the integers $d$ and $\delta$ from
(ii), (iii) are equal; we call this
common value the {\it diameter} of the
pair.
By \cite[Theorem 10.1]{TD00} the pair $A, A^*$ satisfy two polynomial
equations called the tridiagonal relations;
these generalize the $q$-Serre relations
  \cite[Example~3.6]{qSerre}
and the Dolan-Grady relations
  \cite[Example~3.2]{qSerre}.
See
\cite{bas6,
N:aw,
tersub3,
qSerre,aw}
 for results on the tridiagonal relations.
An ordering of the eigenspaces of $A$ (resp. $A^*$)
is said to be {\em standard} whenever it satisfies 
(\ref{eq:t1})
 (resp. (\ref{eq:t2})). 
We comment on the uniqueness of the standard ordering.
Let $\{V_i\}_{i=0}^d$ denote a standard ordering of the eigenspaces of $A$.
By \cite[Lemma~2.4]{TD00}, 
 the ordering $\{V_{d-i}\}_{i=0}^d$ is also standard and no further
 ordering
is standard.
A similar result holds for the eigenspaces of $A^*$.
Let $\{V_i\}_{i=0}^d$ (resp.
$\{V^*_i\}_{i=0}^d$)
denote a standard ordering of the eigenspaces
 of $A$ (resp. $A^*$).
By \cite[Corollary~5.7]{TD00}, 
for $0 \leq i \leq d$ the spaces $V_i$, $V^*_i$
have the same dimension; we denote
this common dimension by $\rho_i$. 
By \cite[Corollaries 5.7, 6.6]{TD00}
the sequence $\{\rho_i\}_{i=0}^d$ is symmetric and unimodal;
that is $\rho_i=\rho_{d-i}$ for $0 \leq i \leq d$ and
$\rho_{i-1} \leq \rho_i$ for $1 \leq i \leq d/2$.
We call the sequence $\{\rho_i\}_{i=0}^d$ the {\em shape}
of $A,A^*$.
The {\it shape conjecture} 
\cite[Conjecture~13.5]{TD00}
states that if 
$\K$ is algebraically closed then
$\rho_{i} \leq 
\binom{d}{i}
$ for $0 \leq i \leq d$.
The shape conjecture has been proven
for a number of special cases
\cite{
NN,
 IT:Krawt,
IT:aug,
shape}.
The TD pair $A,A^*$ is called {\it sharp} whenever
$\rho_0=1$.
By
\cite[Theorem~1.3]{nomstructure},
if $\K$ is algebraically closed then
 $A,A^*$ is sharp.
It is an open problem to classify the
sharp TD pairs up to isomorphism,
but progress is being made
\cite{
 IT:Krawt,
IT:aug,
NN,
Vidar,
nomtowards,
nom:mu}.
The TD pairs of shape
 $(1,1,\ldots,1)$ are called
{\em Leonard pairs}
 \cite[Definition 1.1]{LS99},
and these are classified 
up to isomorphism \cite{LS99,
TLT:array}.
This classification yields a  
 correspondence between the Leonard pairs and a
family of orthogonal polynomials consisting of the $q$-Racah polynomials
and their relatives 
\cite{AWil,
qrac}.
This family coincides with the terminating branch of the Askey scheme
\cite{KoeSwa}.
See 
\cite{NT:balanced,NT:formula,NT:det,NT:mu,
NT:span,NT:switch,
nomsplit,
madrid}
and the references therein for 
results on Leonard pairs.

\medskip
\noindent We now summarize the present paper.
For the above TD pair $A,A^*$ 
let $\lbrace V_i\rbrace_{i=0}^d$
(resp.
 $\lbrace V^*_i\rbrace_{i=0}^d$) denote a
standard ordering of the eigenspaces of
$A$ (resp. $A^*$).
For $0 \leq i \leq d$ let 
$\theta_i$ (resp. $\theta^*_i$)
denote the eigenvalue of $A$ (resp. $A^*$) 
for $V_i$ (resp. $V^*_i$). By
\cite[Theorem~11.1]{TD00}
the expressions
\begin{eqnarray}
\label{eq:bp1} 
\frac{\theta_{i-2}-\theta_{i+1}}{\theta_{i-1}-\theta_i},  \qquad\qquad
  \frac{\theta^*_{i-2}-\theta^*_{i+1}}{\theta^*_{i-1}-\theta^*_i}
\end{eqnarray}
are equal and independent of $i$ for $2 \leq i \leq d-1$.
For this constraint the ``most general'' solution is
\begin{eqnarray}
&&\theta_i = a + bq^{2i-d} + cq^{d-2i}
 \qquad \qquad (0 \leq i \leq d),
\label{eq:cond1}
\\
&&\theta^*_i = a^* + b^*q^{2i-d} + c^*q^{d-2i}
 \qquad \qquad (0 \leq i \leq d),
\label{eq:cond2}
\\
&& q, \;a, \;b, \;c, \;a^*, \;b^*, \;c^* \in {\overline \F},
\label{eq:cond3}
\\
&&q \not=0, \qquad q^2 \not=1, \qquad q^2 \not=-1, \qquad bb^*cc^*\not=0.
\label{eq:cond4}
\end{eqnarray}
For this solution $q^2 + q^{-2}+1$ is the common value of
(\ref{eq:bp1}). The TD pair
 $A,A^*$ is said to have {\it $q$-Racah type} whenever 
(\ref{eq:cond1})--(\ref{eq:cond4}) hold.
 By \cite[Theorem~5.16]{TLT:array}
 the Leonard pairs of $q$-Racah type 
correspond to the 
$q$-Racah polynomials.

\medskip
\noindent
In this paper we classify up to isomorphism
the TD
pairs over an algebraically closed field that have $q$-Racah type.
Our main result is Theorem 
 \ref{thm:mainth}.
The proof involves the representation theory
of the quantum affine algebra
 $U_q(\widehat{ \mathfrak{sl}}_2)$.

\section{Tridiagonal systems}

\indent
When working with a TD pair, it is often convenient to consider
a closely related object called a TD system.
To define a TD system, we recall a few concepts from linear
algebra.
Let $V$ denote a vector space over $\K$ with finite
positive dimension.
Let ${\rm End}(V)$ denote the $\K$-algebra of all linear
transformations from $V$ to $V$.
Let $A$ denote a diagonalizable element of $\mbox{\rm End}(V)$.
Let $\{V_i\}_{i=0}^d$ denote an ordering of the eigenspaces of $A$
and let $\{\theta_i\}_{i=0}^d$ denote the corresponding ordering of
the eigenvalues of $A$.
For $0 \leq i \leq d$ define $E_i \in 
\mbox{\rm End}(V)$ 
such that $(E_i-I)V_i=0$ and $E_iV_j=0$ for $j \neq i$ $(0 \leq j \leq d)$.
Here $I$ denotes the identity of $\mbox{\rm End}(V)$.
We call $E_i$ the {\em primitive idempotent} of $A$ corresponding to $V_i$
(or $\theta_i$).
Observe that
(i) $\sum_{i=0}^d E_i = I$;
(ii) $E_iE_j=\delta_{i,j}E_i$ $(0 \leq i,j \leq d)$;
(iii) $V_i=E_iV$ $(0 \leq i \leq d)$;
(iv) $A=\sum_{i=0}^d \theta_i E_i$.
Moreover
\begin{equation}         \label{eq:defEi}
  E_i=\prod_{\stackrel{0 \leq j \leq d}{j \neq i}}
          \frac{A-\theta_jI}{\theta_i-\theta_j}.
\end{equation}
Note that each of 
$\{A^i\}_{i=0}^d$,
$\{E_i\}_{i=0}^d$ is a basis for the $\K$-subalgebra
of $\mbox{\rm End}(V)$ generated by $A$.
Moreover $\prod_{i=0}^d(A-\theta_iI)=0$.
Now let $A,A^*$ denote a TD pair on $V$.
An ordering of the primitive idempotents 
 of $A$ (resp. $A^*$)
is said to be {\em standard} whenever
the corresponding ordering of the eigenspaces of $A$ (resp. $A^*$)
is standard.

\begin{definition}
{\rm \cite[Definition~2.1]{TD00}}
 \label{def:TDsystem} 
\rm
Let $V$ denote a vector space over $\K$ with finite
positive dimension.
By a {\it tridiagonal system} (or {\it  $TD$ system}) on $V$ we mean a sequence
\[
 \Phi=(A;\{E_i\}_{i=0}^d;A^*;\{E^*_i\}_{i=0}^d)
\]
that satisfies (i)--(iii) below.
\begin{itemize}
\item[(i)]
$A,A^*$ is a TD pair on $V$.
\item[(ii)]
$\{E_i\}_{i=0}^d$ is a standard ordering
of the primitive idempotents of $A$.
\item[(iii)]
$\{E^*_i\}_{i=0}^d$ is a standard ordering
of the primitive idempotents of $A^*$.
\end{itemize}
We say $\Phi$ is {\em over} $\K$.
We call $V$ the {\it underlying vector space}.
\end{definition}

\medskip
\noindent The notion of isomorphism for TD systems
is defined in
\cite[Section~3]{nomsharp}.

\begin{definition}        \label{def}
\rm
Let $\Phi=(A;\{E_i\}_{i=0}^d;A^*$; $\{E^*_i\}_{i=0}^d)$ 
denote a TD system on $V$.
For $0 \leq i \leq d$ let $\theta_i$ (resp. $\theta^*_i$)
denote the eigenvalue of $A$ (resp. $A^*$)
associated with the eigenspace $E_iV$ (resp. $E^*_iV$).
We call $\{\theta_i\}_{i=0}^d$ (resp. $\{\theta^*_i\}_{i=0}^d$)
the {\em eigenvalue sequence}
(resp. {\em dual eigenvalue sequence}) of $\Phi$.
We observe $\{\theta_i\}_{i=0}^d$ (resp. $\{\theta^*_i\}_{i=0}^d$) are
mutually distinct and contained in $\K$.
We say $\Phi$ is
{\it sharp} whenever the TD pair
$A,A^*$ is sharp.
\end{definition}
\medskip

\noindent 
The following notation will be useful.

\begin{definition}
\rm
Let $\lambda$ denote an indeterminate and let $\K[\lambda]$
denote the $\K$-algebra consisting of the polynomials
in $\lambda$ that have all coefficients in $\K$.
Let $\lbrace \theta_i\rbrace_{i=0}^d$  and
 $\lbrace \theta^*_i\rbrace_{i=0}^d$
denote scalars in $\K$.
Then for $0 \leq i \leq d$ we define the following polynomials in
$\K[\lambda]$:
\begin{eqnarray*}
 \tau_i &=& 
  (\lambda-\theta_0)(\lambda-\theta_1)\cdots(\lambda -\theta_{i-1}), \\
 \eta_i &=&
  (\lambda-\theta_d)(\lambda-\theta_{d-1})\cdots(\lambda-\theta_{d-i+1}),  \\
 \tau^*_i &=&
  (\lambda-\theta^*_0)(\lambda-\theta^*_1)\cdots(\lambda-\theta^*_{i-1}), \\
 \eta^*_i &=&
  (\lambda-\theta^*_d)(\lambda-\theta^*_{d-1})\cdots(\lambda-\theta^*_{d-i+1}).
\end{eqnarray*}
Note that each of $\tau_i$, $\eta_i$, $\tau^*_i$, $\eta^*_i$ is
monic with degree $i$.
\end{definition}

\noindent 
We now recall the split sequence of a sharp TD system.
This sequence was originally defined
in \cite[Section~5]{IT:Krawt} using the split
decomposition \cite[Section~4]{TD00},
 but
in \cite{nom:mu}
an alternate
definition was introduced that is more convenient
to our purpose.

\begin{definition}
\label{def:split}
\rm
{\rm \cite[Definition~2.5]{nom:mu}}
Let 
$(A; \lbrace E_i\rbrace_{i=0}^d; A^*; \lbrace E^*_i\rbrace_{i=0}^d)$
denote a sharp TD system over $\K$, with eigenvalue
sequence $\lbrace \theta_i \rbrace_{i=0}^d$
and dual eigenvalue sequence
 $\lbrace \theta^*_i \rbrace_{i=0}^d$.
By 
\cite[Lemma 5.4]{nomstructure},
 for $0 \leq i \leq d$
there exists a unique $\zeta_i \in \K$ such that 
\begin{eqnarray*}
E^*_0 \tau_i(A) E^*_0 = 
\frac{\zeta_i E^*_0}{
(\theta^*_0-\theta^*_1) 
(\theta^*_0-\theta^*_2) 
\cdots
(\theta^*_0-\theta^*_i)}. 
\end{eqnarray*}
Note that $\zeta_0=1$.
We call $\lbrace \zeta_i \rbrace_{i=0}^d$
the {\it split sequence} of the TD system.
\end{definition}

\begin{definition}
\label{def:pa}
\rm
Let $\Phi$
denote a sharp TD system.
By the {\it parameter array} of $\Phi$ we mean the
sequence
 $(\{\theta_i\}_{i=0}^d; \{\theta^*_i\}_{i=0}^d; \{\zeta_i\}_{i=0}^d)$
where
 $\{\theta_i\}_{i=0}^d$
(resp. 
$\{\theta^*_i\}_{i=0}^d$
)
is the eigenvalue sequence 
(resp. dual eigenvalue sequence)
of $\Phi$ and
$\{\zeta_i\}_{i=0}^d$ is the split sequence of $\Phi$.
\end{definition}

\noindent 
The following result shows the significance
of the parameter array. 

\begin{proposition}
\label{thm:isopa}
{\rm 
\cite{IT:aug},
\cite[Theorem~1.6]{nomstructure}}
Two sharp TD systems over $\K$ are isomorphic
if and only if they have the same parameter array.
\end{proposition}

\section{The classification} 

\noindent In  this section we state our main result
and discuss its significance. 

\begin{definition}
\rm
\label{def:qrac}
Let $d$ denote a nonnegative integer and let
 $(\{\theta_i\}_{i=0}^d; \{\theta^*_i\}_{i=0}^d)$
denote a sequence of scalars taken from $\K$. We call
this sequence {\it $q$-Racah} whenever 
the following {\rm (i)}, {\rm (ii)} hold.
\begin{itemize}
\item[\rm (i)]
$\theta_i \neq \theta_j$, $\theta^*_i \neq \theta^*_j$ 
if $i \neq j$ $(0 \leq i,j \leq d)$.
\item[\rm (ii)] There exist $q,a,b,c,a^*,b^*,c^*$ 
that satisfy (\ref{eq:cond1})--(\ref{eq:cond4}).
\end{itemize}
\end{definition}

\noindent Referring to Definition
\ref{def:qrac}, condition (i) implies
a restriction on the scalars in condition (ii). We now
clarify this restriction.

\begin{lemma}
\label{lem:clar}
The following are equivalent for all
integers $d\geq 0$, nonzero $q \in {\overline \F}$,
and $a,b,c \in {\overline \F}$:
\begin{itemize}
\item[\rm (i)]
The scalars $\lbrace
a + bq^{2i-d}+cq^{d-2i}
\rbrace_{i=0}^d$ are mutually distinct;
\item[\rm (ii)]
$q^{2i}\not=1$ for $1 \leq i \leq d$ and 
$b\not=cq^{2d-2i}$ for $1\leq i \leq 2d-1$.
\end{itemize}
\end{lemma}
\noindent {\it Proof:} Routine.
\hfill $\Box$ \\

\noindent The following is our main result.

\begin{theorem}
 \label{thm:mainth} 
Assume the field $\K$ is algebraically 
closed and let 
 $d$ denote a nonnegative integer.
Let 
 $(\{\theta_i\}_{i=0}^d; \{\theta^*_i\}_{i=0}^d)$
denote a $q$-Racah sequence of scalars in $\K$
and let 
$\{\zeta_i\}_{i=0}^d$ denote any sequence of scalars
in $\K$. Then the following are equivalent:
\begin{itemize}
\item[\rm (i)]
There exists a TD system $\Phi$ over $\K$ that
has parameter array
 $(\{\theta_i\}_{i=0}^d; \{\theta^*_i\}_{i=0}^d; \{\zeta_i\}_{i=0}^d)$;
\item[\rm (ii)]
$\zeta_0=1$, $\zeta_d \neq 0$, and
\begin{equation}   
        \label{eq:ineq}
0 \neq \sum_{i=0}^d \eta_{d-i}(\theta_0)\eta^*_{d-i}(\theta^*_0) 
\zeta_i.
\end{equation}
\end{itemize}
Suppose {\rm (i)}, {\rm (ii)} hold. 
Then $\Phi$ is unique up to isomorphism of
TD systems.
\end{theorem}
\noindent Our proof of Theorem
 \ref{thm:mainth} is contained in Section 10.

\medskip
\noindent We now discuss the significance of
Theorem \ref{thm:mainth}. 
The following conjectured classification of TD pairs
was introduced in 
\cite[Conjecture~14.6]{IT:Krawt}; see also
\cite[Conjecture~6.3]{nomsharp} and
\cite[Conjecture~11.1]{nom:mu}.

\begin{conjecture}{\rm \cite[Conjecture~14.6]{IT:Krawt}}  
 \label{conj:main}  
Assume the field $\K$ is algebraically closed.
Let $d$ denote a nonnegative integer and  let
\begin{equation}         \label{eq:parrayc}
 (\{\theta_i\}_{i=0}^d; \{\theta^*_i\}_{i=0}^d; \{\zeta_i\}_{i=0}^d)
\end{equation}
denote a sequence of scalars taken from $\K$.
Then there exists a TD system $\Phi$ over $\K$
with parameter array \eqref{eq:parrayc} if and only if
{\rm (i)}--{\rm (iii)} hold below:
\begin{itemize}
\item[\rm (i)]
$\theta_i \neq \theta_j$, 
$\theta^*_i \neq \theta^*_j$ if $i \neq j$ $(0 \leq i,j \leq d)$.
\item[\rm (ii)]
$\zeta_0=1$, $\zeta_d \neq 0$, and
\begin{equation*}   
        \label{eq:ineqc}
0 \neq \sum_{i=0}^d \eta_{d-i}(\theta_0)\eta^*_{d-i}(\theta^*_0) \zeta_i.
\end{equation*}
\item[\rm (iii)]
The expressions
\begin{equation*} 
\label{eq:betaplusonec} 
\frac{\theta_{i-2}-\theta_{i+1}}{\theta_{i-1}-\theta_i},  \qquad\qquad
  \frac{\theta^*_{i-2}-\theta^*_{i+1}}{\theta^*_{i-1}-\theta^*_i}
\end{equation*}
are equal and independent of $i$ for $2 \leq i \leq d-1$.
\end{itemize}
Suppose {\rm (i)}--{\rm (iii)} hold. Then $\Phi$ is unique up to isomorphism of
TD systems.
\end{conjecture}

\noindent
 The ``only if'' direction of
Conjecture 
\ref{conj:main} was proved in
\cite[Section~8]{nomsharp}. 
The last assertion of 
Conjecture 
\ref{conj:main} 
follows from
Proposition \ref{thm:isopa}.
The ``if'' direction of
 Conjecture 
\ref{conj:main} was proved for
 $d\leq 5$  in \cite[Theorem~11.2, Theorem~12.1]{nom:mu}.
Theorem \ref{thm:mainth} establishes
the ``if''
direction of Conjecture 
\ref{conj:main} for the case in which
 $(\{\theta_i\}_{i=0}^d; \{\theta^*_i\}_{i=0}^d)$ has 
$q$-Racah type.
We remark that our forthcoming paper
\cite{IT:aug}
contains
  a comprehensive treatment of the
TD pairs for which $q$ is not a root of unity,
where $q^2+q^{-2}+1$ is the common value
of 
(\ref{eq:bp1}). 
The treatment establishes
the 
 ``if''
direction of Conjecture 
\ref{conj:main} assuming that restriction on $q$.

\section{An outline of the proof for Theorem
\ref{thm:mainth}}

In the proof of 
Theorem \ref{thm:mainth}
the main part is to demonstrate that (ii) implies (i).
This demonstration will consume most of the paper
from Section 6 through Section 10.
Here we summarize the argument.

\medskip
\noindent 
Assuming $\K$ is algebraically closed,
we fix a $q$-Racah sequence
$(\lbrace \theta_i\rbrace_{i=0}^d;
\lbrace \theta^*_i\rbrace_{i=0}^d)$
 of scalars in $\K$, and 
a sequence $\lbrace \zeta_i\rbrace_{i=0}^d$
of scalars in $\K$ that satisfy condition
(ii) of Theorem
 \ref{thm:mainth}. 
Our goal is to display a TD system 
over $\K$ that has parameter array
$(\lbrace \theta_i\rbrace_{i=0}^d;
\lbrace \theta^*_i\rbrace_{i=0}^d; \lbrace \zeta_i\rbrace_{i=0}^d)$.
To this end
we fix $q,a,b,c,a^*,b^*,c^*$ that satisfy
 (\ref{eq:cond1})--(\ref{eq:cond4}).
Associated with $q$ is the  algebra
$U_q(\widehat{ \mathfrak{sl}}_2)$ over $\K$ with Chevalley
generators 
$e^\pm_i$,
$K^{\pm 1}_i $, 
 $i\in \lbrace 0,1\rbrace$. 
We consider 
the $U_q(\widehat{ \mathfrak{sl}}_2)$-modules $V$ of
the form
\begin{eqnarray*}
\label{eq:standard}
V=V(\alpha_1) \otimes
V(\alpha_2) \otimes
\cdots 
\otimes
V(\alpha_d),
\end{eqnarray*}
where each $V(\alpha_i)$ is a 2-dimensional
evaluation module with evaluation parameter
$\alpha_i$.
The module $V$
 decomposes into a direct sum
of weight spaces $\lbrace U_i\rbrace_{i=0}^d$,
with $K_0-q^{2i-d}1$ and $K_1-q^{d-2i}1$ vanishing
on $U_i$ for $0 \leq i \leq d$.
The dimension of
$U_i$ is
$\binom{d}{i}$ for  $0 \leq i \leq d$.
Using $b,c,b^*,c^*$ we obtain certain elements
$R,L \in 
U_q(\widehat{ \mathfrak{sl}}_2)$
such that
\begin{eqnarray*}
 RU_i \subseteq U_{i+1}, \qquad \qquad 
 LU_i \subseteq U_{i-1}
\qquad \qquad (0 \leq i \leq d),
\end{eqnarray*}
where $U_{-1}=0$ and
$U_{d+1}=0$.
The $R,L$ satisfy 
some attractive equations and act nicely
with respect to the Hopf algebra
structure.
Using $R,L$ we associate with
$V$ a certain monic univariate polynomial
$P_V$ of degree $d$, called
the (nonstandard) Drinfel'd polynomial
\cite{drin}.
Using some properties of 
$P_V$ we show  that 
the parameters $\lbrace \alpha_i\rbrace_{i=1}^d$
can be chosen so that
$\zeta_i$ is the eigenvalue of $L^iR^i$ on 
$U_0$ for $0 \leq i \leq d$; for the rest of
this section we work
with this choice. 
Define
\begin{eqnarray*}
A &=& a 1 + bK_0 + cK_1 + R,
\\
A^* &=& a^* 1 + b^*K_0 + c^*K_1 + L.
\end{eqnarray*}
Using the equations satisfied by $R,L$ we show that
$A,A^*$ satisfy a pair of
tridiagonal relations.
From the construction
\begin{eqnarray*}
(A-\theta_i 1)U_i \subseteq U_{i+1},
\qquad \qquad 
(A^*-\theta^*_i 1)U_i \subseteq U_{i-1},
\qquad \qquad (0 \leq i \leq d).
\end{eqnarray*}
Using this we argue that $A$ (resp. $A^*$) is 
diagonalizable 
on $V$ with eigenvalues $\lbrace \theta_i\rbrace_{i=0}^d$
(resp. $\lbrace \theta^*_i\rbrace_{i=0}^d$).
For $0 \leq i \leq d$ let $E_i$ (resp. $E^*_i$)
denote the element of 
$U_q(\widehat{ \mathfrak{sl}}_2)$
that acts on $V$ as the primitive idempotent
of $A$ (resp. $A^*$) associated with
$\theta_i$ (resp. $\theta^*_i$).
Using the tridiagonal relations
and a few other facts 
we show that
on $V$,
\begin{eqnarray}
\label{eq:triple}
E_iA^*E_j=0,
\qquad \quad 
E^*_iAE^*_j=0
\qquad \quad 
\mbox{\rm if} \quad|i-j|>1,
\qquad \quad
(0 \leq i,j\leq d).
\end{eqnarray}
Using
 $\zeta_d\not=0$ and
        (\ref{eq:ineq})
we show that on $V$,
\begin{eqnarray}
\label{eq:000}
E^*_0E_dE^*_0 \not=0, 
\qquad \qquad
E^*_0E_0E^*_0 \not=0. 
\end{eqnarray}
Let $T$ denote the subalgebra of 
$U_q(\widehat{ \mathfrak{sl}}_2)$
generated by $A,A^*$
and let $TE^*_0V$ denote
the $T$-submodule of $V$ generated by
$E^*_0V=U_0$.
We show that $TE^*_0V$ contains a unique maximal proper $T$-submodule;
denote this by $M$ and consider the quotient
$T$-module $L=TE^*_0V/M$.
By construction the $T$-module
$L$ is irreducible.
Using this and 
(\ref{eq:triple}),
 (\ref{eq:000})
we show that
 the elements
$(A;\lbrace E_i\rbrace_{i=0}^d;
A^*;\lbrace E^*_i\rbrace_{i=0}^d)$
act on $L$
 as a TD system which we denote by $\Phi$.
 By the construction 
$\Phi$ has eigenvalue sequence
$\lbrace \theta_i\rbrace_{i=0}^d$
and dual eigenvalue sequence
$\lbrace \theta^*_i\rbrace_{i=0}^d$.
We argue using the choice of $V$
that $\Phi$ has split
sequence
$\lbrace \zeta_i\rbrace_{i=0}^d$.
Therefore  $\Phi$
 has parameter
array
$(\lbrace \theta_i\rbrace_{i=0}^d;
\lbrace \theta^*_i\rbrace_{i=0}^d; \lbrace \zeta_i\rbrace_{i=0}^d)$
and we have accomplished our goal.

\section{The algebra $U_q(\widehat{ \mathfrak{sl}}_2)$}

\noindent
In this section we recall some facts
about 
$U_q(\widehat{ \mathfrak{sl}}_2)$
that we will use in the proof of
Theorem
\ref{thm:mainth}.
We  follow the notational conventions
of Chari and Pressley
\cite{charp}, \cite{charpbook}.

\medskip
\noindent
Throughout this section
assume $\K$ is algebraically closed.
We fix a nonzero $q \in \K$ such that $q^2\not=1$,
and adopt the
following notation:
\begin{eqnarray}
\lbrack n \rbrack_q = \frac{q^n-q^{-n}}{q-q^{-1}}
\qquad \qquad n = 0,1,2,\ldots 
\label{eq:nbrack}
\end{eqnarray}

\begin{definition}
\label{def:uq}
\rm
\cite[p.~262]{charp}
The quantum affine algebra
$U_q(\widehat{ \mathfrak{sl}}_2)$ is the associative $\K$-algebra
with $1$, defined by
generators $e^{\pm}_i$, $K_i^{{\pm}1}$, $i\in \lbrace 0,1\rbrace $
and the following relations:
\begin{eqnarray}
K_iK^{-1}_i &=& 
K^{-1}_iK_i =  1,
\label{eq:buq1}
\\
K_0K_1&=& K_1K_0,
\label{eq:buq2}
\\
K_ie^{\pm}_iK^{-1}_i &=& q^{{\pm}2}e^{\pm}_i,
\label{eq:buq3}
\\
K_ie^{\pm}_jK^{-1}_i &=& q^{{\mp}2}e^{\pm}_j, \qquad i\not=j,
\label{eq:buq4}
\\
\lbrack e^+_i, e^-_i\rbrack &=& \frac{K_i-K^{-1}_i}{q-q^{-1}},
\label{eq:buq5}
\\
\lbrack e^{\pm}_0, e^{\mp}_1\rbrack &=& 0,
\label{eq:buq6}
\end{eqnarray}
\begin{eqnarray}
(e^{\pm}_i)^3e^{\pm}_j -  
\lbrack 3 \rbrack_q (e^{\pm}_i)^2e^{\pm}_j e^{\pm}_i 
+\lbrack 3 \rbrack_q e^{\pm}_ie^{\pm}_j (e^{\pm}_i)^2 - 
e^{\pm}_j (e^{\pm}_i)^3 =0, \qquad i\not=j.
\label{eq:buq7}
\end{eqnarray}
In 
(\ref{eq:buq5}),
(\ref{eq:buq6}) the expression
$\lbrack r,s\rbrack$ means $rs-sr$.
We call $e^{\pm}_i$, $K_i^{{\pm}1}$, $i\in \lbrace 0,1\rbrace $
the {\it Chevalley generators} for
$U_q({\widehat{\mathfrak{sl}}}_2)$.
\end{definition}

\noindent In \cite[Section~4]{charp}
Chari and Pressley
consider some finite-dimensional
irreducible 
$U_q(\widehat{ \mathfrak{sl}}_2)$-modules
$V_n(\alpha)$, where $\alpha$ is a nonzero scalar in 
$\F$ and $n$ is a positive integer.
These modules are called evaluation modules.
 The scalar $\alpha$ is the evaluation parameter
and $n+1$ is the dimension. We will make
use of $V_1(\alpha)$; for notational convenience
we denote this module by $V(\alpha)$.

\begin{definition}
\label{lem:evmod}
{\rm \cite[Section~4]{charp}}
\rm
For all nonzero $\alpha \in \K$ 
the $U_q(\widehat{ \mathfrak{sl}}_2)$-module $V(\alpha)$
has a basis $x,y$ on which
the Chevalley generators act as follows:
\begin{eqnarray*}
K_1 x&=& qx, \quad \qquad \qquad
K_1 y\;=\; q^{-1}y,
\\
e_1^- x &=& y, \qquad 
\qquad \qquad  e_1^- y \;=\; 0,
\\
e_1^+ x &=& 0, \qquad \qquad \qquad
e_1^+ y 
\;=\;
x,
\\
K_0 x &=&q^{-1}x, \qquad \qquad K_0y \;=\;qy,
\\
e_0^- x &=& 0, \qquad \qquad 
\;\;\quad e_0^- y \;=\; 
 q\alpha^{-1}x,
\\
e_0^+ x &=& q^{-1}\alpha y, 
\quad  \;
\qquad \;e_0^+ y \;=\; 0.
\end{eqnarray*}
\end{definition}

\noindent 
We now recall how the tensor product of
two 
$U_q(\widehat{ \mathfrak{sl}}_2)$-modules becomes
 a 
$U_q(\widehat{ \mathfrak{sl}}_2)$-module.
In what follows all unadorned tensor products are
meant to be over $\F$.

\begin{lemma}
{\rm \cite[p.~263]{charp}}
\label{lem:hopf} 
$U_q({\widehat{\mathfrak{sl}}}_2)$ has the following Hopf algebra
structure. The comultiplication $\Delta$ satisfies
\begin{eqnarray*}
\Delta(e_i^{+})&=& e_i^+\otimes K_i+1 \otimes e_i^+,
\\
\Delta(e_i^{-})&=& e_i^-\otimes 1+ K_i^{-1}\otimes e_i^-,
\\
\Delta(K_i)&=& K_i\otimes K_i.
\end{eqnarray*}
The counit $\varepsilon$ satisfies
\begin{eqnarray*}
\varepsilon(e_i^{\pm})=0, \qquad \qquad
\varepsilon(K_i)=1.
\end{eqnarray*}
The antipode $S$ satifies
\begin{eqnarray*}
S(K_i)= K_i^{-1}, 
\qquad 
\qquad
S(e_i^+)= -e_i^+K_i^{-1}, 
\qquad 
\qquad
S(e_i^-)= -K_ie_i^-.
\end{eqnarray*}
\end{lemma}

\noindent Combining Lemma \ref{lem:hopf}
with 
\cite[p.~110]{charpbook} we routinely obtain the following.

\begin{lemma}
\label{lem:hopfmod}
Let $V,W$ denote 
$U_q(\widehat{ \mathfrak{sl}}_2)$-modules.
Then the tensor product
$V\otimes W$ has  
the following 
$U_q(\widehat{ \mathfrak{sl}}_2)$-module 
structure. 
For $v \in V$, for $w \in W$ and for $i \in \lbrace 0,1\rbrace $, 
\begin{eqnarray*}
e_i^{+}(v\otimes w)&=& e_i^+v\otimes K_iw+v \otimes e_i^+w,
\\
e_i^{-}(v\otimes w)&=& e_i^-v\otimes w+ K_i^{-1}v\otimes e_i^-w,
\\
K_i(v\otimes w)&=& K_iv\otimes K_iw.
\end{eqnarray*}
\end{lemma}

\begin{definition}
\label{def:trivmod}
\rm
\cite[p.~110]{charpbook}
There exists
a one dimensional
$U_q(\widehat{ \mathfrak{sl}}_2)$-module on which
each element $z \in
U_q(\widehat{ \mathfrak{sl}}_2)$
acts as $\varepsilon(z)I$,
where $\varepsilon$ is
from 
Lemma 
\ref{lem:hopf} and $I$ is the identity map.
In particular on this module 
each of
$e^{\pm}_0$,
$e^{\pm}_1$
vanishes
and each of
$K^{\pm 1}_0$, 
$K^{\pm 1}_1$ 
acts as $I$. This module
is irreducible and
 unique up
to isomorphism.
We call this module
the
 {\it trivial}
$U_q(\widehat{ \mathfrak{sl}}_2)$-module.
\end{definition}

\begin{definition}
\rm
\label{def:st}
Let $d$ denote a nonnegative integer.
By a {\it standard 
$U_q(\widehat{ \mathfrak{sl}}_2)$-module
of diameter $d$} we mean 
\begin{eqnarray}
 V(\alpha_1 ) \otimes 
 V(\alpha_2 ) \otimes 
\cdots 
\otimes
V(\alpha_d),
\label{eq:tensor}
\end{eqnarray}
where $0 \not=\alpha_i \in \K$ for $1 \leq i \leq d$.
For $d=0$ we interpret 
(\ref{eq:tensor})
 to be the trivial 
$U_q(\widehat{ \mathfrak{sl}}_2)$-module.
\end{definition}
For our purpose we only need those
standard
$U_q(\widehat{ \mathfrak{sl}}_2)$-modules
 of diameter $d$ such that $q^{2i}\not=1 $ for
$1 \leq i \leq d$. The following definition
will facilitate our discussion of these modules.
\begin{definition}
\label{def:feas}
\rm
An integer $d$ will be called {\it feasible} (with respect to $q$)
 whenever
$d\geq 0$ and $q^{2i}\not=1$ for $1 \leq i \leq d$.
Note that $0$ and $1$ are feasible.
\end{definition}

\noindent The following result is immediate from
Definition
\ref{def:feas}.

\begin{lemma}
\label{lem:mutdist}
For all feasible integers $d$ the scalars 
$\lbrace q^{d-2i}\rbrace_{i=0}^d$ are mutually distinct.
\end{lemma}

\noindent Let $V$ denote a standard 
$U_q(\widehat{ \mathfrak{sl}}_2)$-module
with feasible diameter $d$.
The $\K$-vector space $V$ has a basis 
\begin{eqnarray}
v_1 \otimes v_2 \otimes \cdots \otimes v_d
\qquad \qquad v_i \in \lbrace x,y\rbrace \qquad \qquad (1 \leq i \leq d).
\label{eq:2dbasis}
\end{eqnarray}
Note that $V$ has dimension $2^d$.
For notational convenience we abbreviate the basis
(\ref{eq:2dbasis}) 
 as
follows.
For  all subsets $s \subseteq 
\lbrace 1,2,\ldots, d\rbrace$
define $u_s=
v_1 \otimes v_2 \otimes \cdots \otimes v_d$,
where $v_i=x$ if $i \not\in s$ and 
$v_i=y$ if $i \in s$ $(1 \leq i \leq d)$.
For example $u_\emptyset=x \otimes x \otimes \cdots \otimes x$.
Pick a subset 
 $s \subseteq 
\lbrace 1,2,\ldots, d\rbrace $.
By Lemma \ref{lem:hopfmod} we have
\begin{eqnarray}
K_0 u_s = 
q^{2|s|-d}u_s,
\qquad \qquad 
K_1 u_s = 
q^{d-2|s|}u_s,
\label{eq:kus}
\end{eqnarray}
 where $|s|$ denotes the
cardinality of $s$.
Thus each of $K_0$, $K_1$ is diagonalizable on $V$ 
 with eigenvalues
$\lbrace q^{d-2i}\rbrace_{i=0}^d$.
Moreover 
$K_0, K_1$ are inverses of one another on $V$.
For $0 \leq i \leq d$ define
\begin{eqnarray*}
U_i = \mbox{\rm Span}\lbrace u_s \;|\;s \subseteq 
\lbrace 1,2,\ldots, d\rbrace, \quad |s|=i\rbrace.
\end{eqnarray*}
Note that 
$V =\sum_{i=0}^d U_i$ (direct sum), and that  
$U_i$ has dimension  
$\binom{d}{i}$ for  $0 \leq i \leq d$.
Moreover 
\begin{eqnarray}
\label{eq:nnkmove}
&&(K_0-q^{2i-d}1)U_i=0,
\qquad \qquad  
 (K_1-q^{d-2i}1)U_i=0
\end{eqnarray}
for $0 \leq i \leq d$.
Combining (\ref{eq:buq3}), 
(\ref{eq:buq4})
with
 (\ref{eq:nnkmove})
we find that for $0 \leq i \leq d$,
\begin{eqnarray}
&&e^+_0 U_i \subseteq U_{i+1}, \qquad \qquad  
e^-_1 U_i \subseteq U_{i+1}, 
\label{eq:altemove1}
\\
&&e^-_0 U_i \subseteq U_{i-1}, \qquad \qquad  
e^+_1 U_i \subseteq U_{i-1}, 
 \label{eq:altemove2}
\end{eqnarray}
where $U_{-1}=0$ and $U_{d+1}=0$.
We call the sequence
$\lbrace U_i\rbrace_{i=0}^d$ the
{\it weight space decomposition} of $V$.
We call $u_\emptyset$ the {\it highest weight vector} of $V$.
The action
of $e^\pm_0$, 
$e^\pm_1$ on the basis
(\ref{eq:2dbasis}) can be obtained using
 Lemma
\ref{lem:hopfmod} but the answer is slightly complicated.
For all subsets 
$s \subseteq \lbrace 1,\ldots, d\rbrace$
let
$\overline s$ 
denote the complement of
$s$ in 
$\lbrace 1,\ldots, d\rbrace$.
By Lemma
\ref{lem:hopfmod},
\begin{eqnarray}
e^-_1u_s&=& \sum_{i \in {\overline s}}
u_{s \cup i} 
q^{
\vert \lbrace 1,2,\ldots, i-1 \rbrace \cap s\vert
-
\vert \lbrace 1,2,\ldots, i-1 \rbrace \cap {\overline s} \vert
},
\label{eq:e1maction}
\\
e^+_1u_s&=& \sum_{i \in s}
u_{s \backslash i} 
q^{
\vert \lbrace i+1,i+2,\ldots, d \rbrace \cap {\overline s} \vert
-
\vert \lbrace i+1,i+2,\ldots, d \rbrace \cap s \vert
},
\label{eq:e1paction}
\\
e^-_0u_s&=& \sum_{i \in s}
u_{s \backslash i} 
\alpha^{-1}_i
q^{
\vert \lbrace 1,2,\ldots, i-1 \rbrace \cap {\overline s}\vert
-
\vert \lbrace 1,2,\ldots, i-1 \rbrace \cap s \vert +1
},
\label{eq:e0maction}
\\
e^+_0u_s&=& \sum_{i \in {\overline s}}
u_{s \cup i} 
\alpha_i 
q^{
\vert \lbrace i+1,i+2,\ldots, d \rbrace \cap s \vert
-
\vert \lbrace i+1,i+2,\ldots, d \rbrace \cap {\overline s} \vert
-1}.
\label{eq:e0paction}
\end{eqnarray}

\section{The elements $R, L$ of 
$U_q(\widehat{\mathfrak{sl}}\sb 2)$}

\noindent  
From now until the end of Section 7 
we adopt the following assumption.
\begin{assumption}
\label{lem:eigdata}
\rm
We assume the field $\K$ is algebraically closed.
We fix nonzero scalars   
$q, b,c,b^*,c^*$ in $\K$ such that $q^2 \not=1$.
\end{assumption}

\noindent
In this section we define
the elements
$R,L$ of  
$U_q(\widehat{\mathfrak{sl}}\sb 2)$ and discuss
their basic properties.

\begin{definition}
\label{def:RL}
\label{def:uv}
\rm
We define 
\begin{eqnarray}
\label{eq:rdef}
R &=& ue^+_0 + v e^-_1 K_1,
\\
\label{eq:ldef}
L &=& u^*e^+_1 + v^* e^-_0 K_0,
\end{eqnarray}
where $u,v,u^*,v^*$ are any scalars in $\F$ such that
\begin{eqnarray}
uv^*&=& -b b^*q^{-1}(q-q^{-1})^2,
\label{eq:uvs}
\\
vu^*&=& -c c^*q^{-1}(q-q^{-1})^2.
\label{eq:vus}
\end{eqnarray}
Note that $u,v,u^*,v^*$ are nonzero.
\end{definition}

\begin{note} \rm
Referring to 
Definition \ref{def:RL},
the choice of $u,v,u^*,v^*$
is immaterial and we could fix
specific values for the duration of
the paper. But doing so
tends to
obscure the essential relationships
(\ref{eq:uvs}),
(\ref{eq:vus}).
\end{note}

\begin{lemma}
\label{lem:RLK}
We have
\begin{eqnarray}
\label{eq:RLK}
&&K_0 R K^{-1}_0 = q^2 R,
\qquad \qquad K_1 R K^{-1}_1 = q^{-2} R,
\\
\label{eq2:RLK}
&&K_0 L K^{-1}_0 = q^{-2} L,
\qquad \qquad K_1 L K^{-1}_1 = q^{2} L.
\end{eqnarray}
\end{lemma}
\noindent {\it Proof:} 
Use 
(\ref{eq:buq1})--(\ref{eq:buq4}) 
and 
Definition \ref{def:RL}.
\hfill $\Box$ \\

\noindent Let 
$U_q(L(\mathfrak{sl}\sb 2))$
denote
the quotient 
of $U_q(\widehat{\mathfrak{sl}}\sb 2)$
by the two sided ideal
 generated by $K_0K_1-1$.
The name
 $U_q(L(\mathfrak{sl}\sb 2))$
is motivated 
 by the fact 
that this algebra
is a $q$-deformation of the universal enveloping
algebra of the loop algebra
$L(\mathfrak{sl}\sb 2)= 
\mathfrak{sl}\sb 2 \otimes \lbrack t,t^{-1}\rbrack$. 
This is discussed in \cite[Section~3.3]{charp}.
In what follows, we will use the same notation for
an element of 
$U_q(\widehat{\mathfrak{sl}}\sb 2)$
and its image in
$U_q(L(\mathfrak{sl}\sb 2))$.

\begin{lemma}
\label{lem:rlk}
The following hold in 
$U_q(L(\mathfrak{sl}\sb 2))$:
\begin{equation}
\begin{split}
&R^3L-\lbrack 3\rbrack_q R^2LR+
\lbrack 3\rbrack_q RLR^2-LR^3
\\
&\qquad \qquad = \;
(q-q^{-1})
(q^2-q^{-2})
(q^3-q^{-3})
(cc^*K_1R^2K_1 
 - bb^*K_0R^2K_0),
\label{eq:qsrl1}
\end{split}
\end{equation}
\begin{equation}
\begin{split}
&L^3R-\lbrack 3\rbrack_q L^2RL+
\lbrack 3\rbrack_q LRL^2-RL^3
\\
&\qquad \qquad = \; (q-q^{-1})
(q^2-q^{-2})
(q^3-q^{-3})
 (bb^*K_0L^2K_0 - cc^*K_1L^2K_1).
\label{eq:qsrl2}
\end{split}
\end{equation}
\end{lemma}
\noindent {\it Proof:} 
To verify 
(\ref{eq:qsrl1}), 
eliminate $R,L$ using
Definition \ref{def:RL}
and simplify the result using
the relations 
in Definition \ref{def:uq}, together with
the fact in 
$U_q(L(\mathfrak{sl}\sb 2))$ the elements
 $K_0, K_1$ become inverses of one another.
Equation 
(\ref{eq:qsrl2}) is similarly verified.
\hfill $\Box$ \\

\begin{lemma}
\label{lem:rlcom}
For all integers $n\geq 0$ the element $L^nR^n$ commutes
with each of $K_0$, $K_1$.
\end{lemma}
\noindent {\it Proof:} 
Use Lemma
\ref{lem:RLK}.
\hfill $\Box$ \\

\begin{lemma}
\label{lem:rlmove}
Let $V$ denote a standard 
$U_q(\widehat{\mathfrak{sl}}\sb 2)$-module with feasible diameter $d$
and let
$\lbrace U_i\rbrace_{i=0}^d$ denote the 
corresponding weight space decomposition. Then 
\begin{eqnarray}
\label{eq:rlmove}
R U_i \subseteq  U_{i+1},
\qquad \qquad
L U_i \subseteq  U_{i-1} 
\qquad \qquad (0 \leq i \leq d).
\end{eqnarray}
\end{lemma}
\noindent {\it Proof:} 
Use
(\ref{eq:nnkmove})--(\ref{eq:altemove2})
and
Definition
\ref{def:RL}.
\hfill $\Box$ \\

\section{The split sequence of a standard 
$U_q(\widehat{\mathfrak{sl}}\sb 2)$-module}

\noindent Throughout this section
Assumption
\ref{lem:eigdata} remains in effect.
We fix elements $R,L \in
U_q(\widehat{\mathfrak{sl}}\sb 2)$
as in Definition
\ref{def:RL}.

\begin{definition}
\rm
\label{def:splituq}
Let $V$ denote a standard
$U_q(\widehat{\mathfrak{sl}}\sb 2)$-module with feasible diameter $d$
and let
$\lbrace U_i\rbrace_{i=0}^d$ denote the 
corresponding weight space decomposition. 
By Lemma
\ref{lem:rlmove}, for $0 \leq i \leq d$
the space $U_0$ is invariant under $L^iR^i$;
let $\zeta_i$ denote the corresponding
eigenvalue. Note that $\zeta_0=1$.
We call the sequence $\lbrace \zeta_i\rbrace_{i=0}^d$
the {\it split sequence of $V$}.
\end{definition}

\medskip
\noindent 
 Our goal in this section is to
obtain the following result.

\begin{proposition}
\label{cor:splitany}
Let $d$ denote a feasible integer and 
let  $\lbrace \zeta_i\rbrace_{i=0}^d$ denote a sequence of
scalars in $\K$ such that  $\zeta_0=1$.
Then there
exists a standard 
$U\sb q(\widehat{\mathfrak{sl}}\sb 2)$-module
$V$ of diameter $d$ that has split sequence
$\lbrace \zeta_i\rbrace_{i=0}^d$.
\end{proposition}

\noindent 
In order to prove Proposition  
\ref{cor:splitany} we will consider a generating function
involving the split sequence called the
{\it (nonstandard) Drinfel'd
polynomial}. We will define this polynomial shortly.

\begin{definition}
\rm
\label{def:splituq2}
Let $V$ denote a standard
$U_q(\widehat{\mathfrak{sl}}\sb 2)$-module with
feasible diameter $d$.
For $0 \leq i \leq d$ define
\begin{eqnarray}
\sigma_i = \frac{\zeta_i}{
(q-q^{-1})^2 
(q^2-q^{-2})^2 
\cdots 
(q^i-q^{-i})^2},
\label{eq:splitnorm}
\end{eqnarray}
where $\lbrace \zeta_i\rbrace_{i=0}^d$ is the
split sequence of
$V$. The denominator in
(\ref{eq:splitnorm}) is nonzero by
Definition
\ref{def:feas}.
Observe that $\sigma_0=1$.
We call $\lbrace \sigma_i\rbrace_{i=0}^d$ the {\it normalized
split sequence} of $V$.
\end{definition}

\begin{definition}
\label{def:pi}
\rm
For all integers $i\geq 0$
 define $f_i \in \K\lbrack \lambda \rbrack$
by
\begin{eqnarray*}
f_i = b b^*q^{-2i}+c c^*q^{2i}-\lambda,
\end{eqnarray*}
where $q,b,b^*,c,c^*$ are from
Assumption \ref{lem:eigdata}.
\end{definition}

\begin{definition}
\label{def:P}
\rm
Let $V$ denote a standard
$U\sb q(\widehat{\mathfrak{sl}}\sb 2)$-module with feasible 
diameter $d$.
We define a polynomial
 $P_V \in \K\lbrack \lambda \rbrack$
by 
\begin{eqnarray}
P_V =
(-1)^d \sum_{i=0}^d 
 \sigma_{d-i} f_0 f_1 \cdots f_{i-1},
\label{eq:drinp}
\end{eqnarray}
where $\lbrace \sigma_i\rbrace_{i=0}^d$ is the normalized
split sequence of $V$.
We observe that $P_V$ is monic with degree $d$.
We call $P_V$ the {\it (nonstandard) Drinfel'd polynomial}
of $V$. 
\end{definition}

\begin{note} \rm
The Drinfel'd polynomial from
\cite[Definition~9.3]{drin} is equal
to the polynomial $P_V$ from
Definition \ref{def:P}, times $(-1)^d(q-q^{-1})^2(q^2-q^{-2})^2
\cdots (q^d-q^{-d})^2$.
\end{note}

\noindent From now on, when we refer to the Drinfel'd polynomial
we mean the nonstandard Drinfel'd polynomial from
Definition \ref{def:P}.

\medskip
\noindent We now compute the Drinfel'd polynomial for
a few easy cases.

\begin{lemma}
\label{lem:trivP}
Let $V$ denote the trivial 
$U\sb q(\widehat{\mathfrak{sl}}\sb 2)$-module
from Definition
\ref{def:trivmod}. Then $P_V=1$.
\end{lemma}
\noindent {\it Proof:} 
Routine.
\hfill $\Box$ \\

\begin{lemma}
\label{lem2:deq1}
Pick a nonzero $\alpha \in \K$ and consider the
$U\sb q(\widehat{\mathfrak{sl}}\sb 2)$-module
$V=V(\alpha)$ from
Definition
\ref{lem:evmod}. The corresponding 
Drinfel'd polynomial satisfies
\begin{eqnarray}
P_V = \lambda - \frac{\alpha u u^* q^{-2} + \alpha^{-1} v v^* q^2}
{q^{-1}(q-q^{-1})^2},
\label{eq2:pv}
\end{eqnarray}
where $u,v,u^*,v^*$ are from
Definition
\ref{def:RL}.
\end{lemma}
\noindent {\it Proof:}
Let $\sigma_1$ denote term one of the
normalized split sequence for $V$. We show
\begin{eqnarray}
\sigma_1 = \frac{\alpha u u^* q^{-2} + \alpha^{-1} v v^* q^2}
{q^{-1}(q-q^{-1})^2}
- bb^*  - cc^*.
\label{eq2:sigma1}
\end{eqnarray}
By Definitions
\ref{def:splituq} and
\ref{def:splituq2},
$\sigma_1 = \zeta_1(q-q^{-1})^{-2}$ where $\zeta_1$
is the eigenvalue of $LR$ on the weight space $U_0$.
Let $x,y$ denote the basis for 
$V$ from
Definition
\ref{lem:evmod}.
By construction $x$ is a basis for $U_0$,
so $x$ is an eigenvector for $LR$
with eigenvalue $\zeta_1$.
 Using
Definition
\ref{lem:evmod}
 and 
Definition \ref{def:RL} we 
find
$Rx = (uq^{-1}\alpha + vq)y$
and 
$Ly = (u^* + v^*q^2\alpha^{-1})x$;
therefore $\zeta_1 = 
(uq^{-1}\alpha + vq)
(u^* + v^*q^2\alpha^{-1})$.
Evaluating this using
(\ref{eq:uvs}), 
(\ref{eq:vus})
and 
$\sigma_1 = \zeta_1(q-q^{-1})^{-2}$
we  obtain 
(\ref{eq2:sigma1}).
Setting $d=1$ and 
$\sigma_0=1$
in 
(\ref{eq:drinp})
we find
$P_V = -\sigma_1 - f_0$.
Evaluating this using
(\ref{eq2:sigma1}) and $f_0= bb^*+cc^*-\lambda$
we obtain
(\ref{eq2:pv}).
\hfill $\Box$ \\

\noindent Consider a standard 
$U\sb q(\widehat{\mathfrak{sl}}\sb 2)$-module
$V=\otimes_{i=1}^d V(\alpha_i)$
with feasible diameter $d$.
We will now show that the Drinfel'd polynomial
of $V$ is equal to the product of
the Drinfel'd polynomials $\prod_{i=1}^d P_{V(\alpha_i)}$.

\begin{lemma}
\label{lem:coprod}
The comultiplication $\Delta$
from
Lemma \ref{lem:hopf} acts on the elements
$R,L$ as follows.

\begin{eqnarray}
\label{eq:deltaR}
\Delta(R) &=& 1 \otimes R + u e^+_0 \otimes K_0
+ v e^-_1 K_1 \otimes K_1,
\\
\label{eq:deltaL}
\Delta(L) &=& 1 \otimes L + u^* e^+_1 \otimes K_1
+ v^* e^-_0 K_0 \otimes K_0.
\end{eqnarray}
\end{lemma}
\noindent {\it Proof:} 
Use
Lemma \ref{lem:hopf}, 
Definition
\ref{def:RL}, and the fact that $\Delta$ is
an algebra homomorphism.
\hfill $\Box$ \\

\begin{lemma}
\label{lem:deltan}
Let $V$ denote a standard
$U\sb q(\widehat{\mathfrak{sl}}\sb 2)$-module
with diameter $1$.
Let $W$ 
denote any
$U\sb q(\widehat{\mathfrak{sl}}\sb 2)$-module.
Then for all integers $n\geq 1$ the following 
 {\rm (i)}, {\rm (ii)} hold
 on $V\otimes W$:
\begin{itemize}
\item[\rm (i)]
$\Delta(R^n)=1 \otimes R^n + \lbrack  n \rbrack_q R_n$ where 
\begin{eqnarray}
R_n = u  q^{n-1} e^+_0 \otimes R^{n-1}K_0 
\;+\;
v  q^{1-n} e^-_1 K_1\otimes R^{n-1}K_1.
\label{eq:uvpart}
\end{eqnarray}
\item[\rm (ii)]
$\Delta(L^n) = 1\otimes L^n + 
\lbrack n \rbrack_q L_n$
where
\begin{eqnarray}
L_n = u^*   q^{1-n} e^+_1 \otimes K_1 L^{n-1} 
\;+\;
v^* q^{n-1} e^-_0 K_0\otimes K_0 L^{n-1}.
\label{eq:usvspart}
\end{eqnarray}
\end{itemize}
\end{lemma}
\noindent {\it Proof:} (i) The proof is by induction on 
$n$. First assume $n=1$. Then the result is immediate from
(\ref{eq:deltaR}). Next assume $n\geq 2$.
By
(\ref{eq:deltaR}) 
 and since
$\Delta(R^n)= \Delta(R^{n-1})\Delta(R)$,
the expression $\Delta(R^n)-1\otimes R^n$ is equal to
\begin{eqnarray}
\label{eq:1}
\bigl(\Delta(R^{n-1})-1\otimes R^{n-1}\bigr)(1 \otimes R)
\end{eqnarray}
plus
 $u$ times
\begin{eqnarray}
 \Delta(R^{n-1})(e^+_0 \otimes K_0)
\label{eq:2}
\end{eqnarray} 
plus
$v$ times
\begin{eqnarray}
 \Delta(R^{n-1})(e^-_1 K_1\otimes K_1).
\label{eq:3}
\end{eqnarray}
We now find the action of
(\ref{eq:1})--(\ref{eq:3}) on $V \otimes W$.
By induction, on 
$V \otimes W$ the expression 
$\Delta(R^{n-1}) - 
1 \otimes R^{n-1}$ is equal to 
$\lbrack n-1 \rbrack_q 
$
times 
\begin{eqnarray}
u   q^{n-2} e^+_0 \otimes R^{n-2}K_0 
\;+\;
v   q^{2-n} e^-_1 K_1\otimes R^{n-2}K_1.
\label{eq:ind}
\end{eqnarray}
By 
this
 and Lemma
\ref{lem:RLK},
on 
$V \otimes W$ the expression
(\ref{eq:1}) is equal to 
$ \lbrack n-1 \rbrack_q$
times
\begin{eqnarray*}
u  q^n e^+_0 \otimes R^{n-1}K_0 
\;+\;
v   q^{-n} e^-_1 K_1\otimes R^{n-1}K_1.
\label{eq:sum1}
\end{eqnarray*}
By 
(\ref{eq:altemove1})
 and since $V$ has diameter 1, the
elements $(e^+_0)^2$ and $e^-_1e^+_0$ are zero on $V$.
Therefore (\ref{eq:ind}) times 
 $e^+_0 \otimes K_0$ is zero on $V\otimes W$, 
so
 (\ref{eq:2}) is equal to
 $e^+_0 \otimes R^{n-1}K_0$ on $V\otimes W$.
Similarly
(\ref{eq:3}) is equal to
 $e^-_1 K_1 \otimes R^{n-1}K_1$ on $V\otimes W$.
By these comments we routinely obtain 
$\Delta(R^n)=1 \otimes R^n + \lbrack  n \rbrack_q R_n$.
\\
\medskip
\noindent (ii) Similar to the proof of (i) above.
\hfill $\Box $ \\

\begin{proposition}
\label{prop:sigmauw}
Pick a feasible integer $d \geq 1$.
Let $V$ denote a standard
$U\sb q(\widehat{\mathfrak{sl}}\sb 2)$-module
with diameter $1$,
and let $W$ 
denote a 
standard
$U\sb q(\widehat{\mathfrak{sl}}\sb 2)$-module
with diameter $d-1$.
Note that 
the $U\sb q(\widehat{\mathfrak{sl}}\sb 2)$-module
$V\otimes W$ is standard with
diameter $d$.
The normalized split sequence for
$V\otimes W$ is described as follows.
\begin{eqnarray*}
\sigma_0(V\otimes W)&=&1,
\\
\sigma_n(V\otimes W)&=& 
(q^{d-n}-q^{n-d})
(bb^*q^{n-d}-cc^*q^{d-n})\sigma_{n-1}(W)
\\
&& \qquad \qquad +
\quad
\sigma_n(W)\;+\;\sigma_1(V)\sigma_{n-1}(W) 
\qquad \qquad (1 \leq n \leq d-1),
\\
\sigma_d(V\otimes W) &=& \sigma_1(V)\sigma_{d-1}(W).
\end{eqnarray*}
\end{proposition}
\noindent {\it Proof:} 
Let $x$ denote   
the highest weight vector in $V$, and observe
\begin{eqnarray}
\label{eq:etadata}
(K_0-q^{-1}1)x=0, 
\qquad \qquad 
(K_1-q 1)x=0.
\end{eqnarray}
Let $\xi$ denote 
the highest weight
vector for $W$, and observe
\begin{eqnarray}
\label{eq:xidata}
(K_0-q^{1-d} 1)\xi=0, 
\qquad \qquad 
(K_1-q^{d-1} 1)\xi=0.
\end{eqnarray}
 Note that
$x \otimes \xi$ is the highest weight vector for
$V \otimes W$.
We claim that for all integers $n\geq 1$,
\begin{equation}
\label{eq:lnrnpre}
\begin{split}
&L^nR^n(x \otimes \xi)- x \otimes L^nR^n \xi  - 
\lbrack n \rbrack^2_q LR x \otimes L^{n-1}R^{n-1}\xi
\\
&\qquad = \quad (q^n-q^{-n})^2(q^{d-n}-q^{n-d})(bb^*q^{n-d}-cc^*q^{d-n}) 
x \otimes L^{n-1}R^{n-1} \xi.
\end{split}
\end{equation}
To prove the claim we evaluate
the left hand side of 
(\ref{eq:lnrnpre}). The term
 $L^nR^n(x \otimes \xi)$ coincides with
the image of $\Delta(L^n)\Delta(R^n)$ on $ x \otimes \xi$.
Computing this image using Lemma
\ref{lem:deltan} one encounters the terms
$L_n(1 \otimes R^n) $ and
$(1 \otimes L^n)R_n$.
Using
(\ref{eq:etadata}) and 
$e^+_1 x=0$,
$e^-_0 x=0$ we find
$L_n(1 \otimes R^n)$ is zero on $x \otimes \xi$.
Using
(\ref{eq:xidata}) and $L^nR^{n-1}\xi=0$ we find
$(1 \otimes L^n)R_n$ is zero on 
 $x \otimes \xi$.
By these comments and Lemma
\ref{lem:deltan} 
the left hand side of 
(\ref{eq:lnrnpre}) is equal to 
$\lbrack n \rbrack^2_q$ times
\begin{eqnarray}
\label{eq:halfway}
L_n R_n (x \otimes \xi)- LRx \otimes L^{n-1}R^{n-1} \xi.
\end{eqnarray}
The two terms in 
(\ref{eq:halfway}) are evaluated as follows.
To evaluate the first term 
use (\ref{eq:uvpart}) and
(\ref{eq:usvspart}).
To evaluate the second term
expand 
$LRx$ using
Definition \ref{def:RL}.
Now reduce further using
(\ref{eq:buq3}),
(\ref{eq:etadata}),
(\ref{eq:xidata}) along with
Definition
\ref{def:uv}, 
Lemma
\ref{lem:rlcom}
and
\begin{eqnarray*}
(e^-_0 e^+_0-1)x = 0,
\qquad \qquad 
(e^+_1 e^-_1-1)x = 0.
\end{eqnarray*}
The reduction shows that
the left hand side of 
(\ref{eq:lnrnpre}) is equal to the right hand
side of 
(\ref{eq:lnrnpre}).
The  claim is now proved.
As we examine the terms in 
(\ref{eq:lnrnpre}), we note
the following 
from Definitions
\ref{def:splituq},
\ref{def:splituq2}.
For $0 \leq n \leq d$ the vector 
$x \otimes \xi$
is an eigenvector for
$L^n R^n$ with eigenvalue 
\begin{eqnarray*}
(q-q^{-1})^2 (q^2-q^{-2})^2\cdots 
(q^n-q^{-n})^2 \sigma_n(V \otimes W).
\end{eqnarray*}
The vector $x$ is an eigenvector
for 
$LR$ with eigenvalue $(q-q^{-1})^2 \sigma_1(V)$.
For $0 \leq n\leq d-1$ the vector 
$\xi$ is an eigenvector for $L^n R^n$ 
with eigenvalue 
\begin{eqnarray*}
(q-q^{-1})^2 (q^2-q^{-2})^2\cdots 
(q^n-q^{-n})^2 \sigma_n(W).
\end{eqnarray*}
Also $L^dR^d \xi=0$ by Lemma
\ref{lem:rlmove}
 and since $W$ has diameter $d-1$.
Evaluating
(\ref{eq:lnrnpre}) using these facts
we obtain the result.
\hfill $\Box$ \\

\begin{proposition}
\label{prop:puw}
With the notation and assumptions
of 
Proposition
\ref{prop:sigmauw},
we have $P_{V \otimes W} = P_V P_W$.
\end{proposition}
\noindent {\it Proof:}
By Definition
\ref{def:P} we have $P_V=-\sigma_1(V)-f_0$.
Again using
Definition
\ref{def:P},
\begin{eqnarray}
\label{eq:pw}
P_W &=&
(-1)^{d-1} \sum_{j=0}^{d-1} 
 \sigma_{d-1-j}(W) f_0 f_1 \cdots f_{j-1},
\\
P_{V\otimes W}&=&
(-1)^d \sum_{i=0}^d 
 \sigma_{d-i}(V \otimes W) f_0 f_1 \cdots f_{i-1}.
\label{eq:pv}
\end{eqnarray}
Using Definition
\ref{def:pi},
\begin{eqnarray}
f_0 = f_j + (q^j-q^{-j})(bb^*q^{-j}-cc^*q^j)
\qquad \quad  (0 \leq j \leq d-1).
\label{eq:f0fj}
\end{eqnarray}
In equation
(\ref{eq:pw}) 
we multiply both sides by
$P_V$ and use
(\ref{eq:f0fj}) to get
\begin{equation}
\label{eq:puwfirst}
\begin{split}
P_V P_W &= (-1)^d
\sum_{j=0}^{d-1} 
\sigma_{d-1-j}(W)  \\
& \qquad \times \quad  f_0 f_1 \cdots f_{j-1}
\bigl(f_j +
(q^{j}-q^{-j})(bb^*q^{-j}-cc^*q^{j})+\sigma_1(V) \bigr).
\end{split}
\end{equation}
In
(\ref{eq:puwfirst}) 
 the sum is a linear
combination of
 $\lbrace f_0f_1 \cdots f_{i-1}\rbrace_{i=0}^{d}$.
In this linear combination, for $0 \leq i \leq d$ let
$\gamma_i$ denote the coefficient of
 $f_0f_1 \cdots f_{i-1}$. We show
\begin{eqnarray}
\label{eq2:gammap}
\gamma_i = \sigma_{d-i}(V\otimes W).
\end{eqnarray}
First assume $i=0$. Then 
(\ref{eq2:gammap}) holds since both sides equal
$\sigma_1(V)\sigma_{d-1}(W)$.
Next assume $1 \leq i \leq d-1$.
By construction
\begin{eqnarray*}
\label{eq2:ppp}
\gamma_i &=&
(q^{i}-q^{-i})(bb^*q^{-i}-cc^*q^{i})\sigma_{d-1-i}(W)
+\sigma_1(V)\sigma_{d-1-i}(W)
+ \sigma_{d-i}(W).
\end{eqnarray*}
Evaluating this using
Proposition
\ref{prop:sigmauw}
we routinely obtain
(\ref{eq2:gammap}).
Next assume $i=d$. Then
(\ref{eq2:gammap})  holds since both sides equal $1$.
We have verified 
(\ref{eq2:gammap}) for $0 \leq i \leq d$.
Therefore 
\begin{eqnarray*}
P_V P_W = 
(-1)^d \sum_{i=0}^d 
 \sigma_{d-i}(V \otimes W) f_0 f_1 \cdots f_{i-1}.
\end{eqnarray*}
Comparing this  with 
(\ref{eq:pv})
 we obtain $P_{V \otimes W}=P_V P_W$.
\hfill $\Box$ \\

\begin{proposition}
\label{thm2:pvpwmain}
Let $V
=\otimes_{i=1}^d V(\alpha_i)$
 denote a standard 
$U\sb q(\widehat{\mathfrak{sl}}\sb 2)$-module 
with feasible diameter $d$.
Then the Drinfel'd polynomial $P_V$ is given by
\begin{eqnarray}
P_V = \prod_{i=1}^d P_{V(\alpha_i)}.
\label{eq:ans}
\end{eqnarray}
\end{proposition}
\noindent {\it Proof:}
Use Proposition
\ref{prop:puw} and induction on $d$.
\hfill $\Box$ \\

\noindent We will make a few more comments on the Drinfel'd polynomial
and then prove Proposition
\ref{cor:splitany}.

\begin{lemma}
\label{lem2:allp}
For all $r \in \K$ there exists a nonzero
$\alpha \in \K$ such that $\lambda -r$ is the 
Drinfel'd polynomial for $V(\alpha)$.
\end{lemma}
\noindent {\it Proof:} 
 Since $ \K$ is algebraically closed 
and $uu^*vv^*\not=0$, 
there exists a nonzero 
$\alpha \in  \K$ such that the fraction on the right
 in (\ref{eq2:pv}) is equal to $r$. 
 The result follows in view of Lemma \ref{lem2:deq1}.  \hfill $\Box$ \\

\begin{note} \rm Referring to Lemma
\ref{lem2:allp}, for a given $r$ the scalar $\alpha$ is not uniquely
determined in general. If $\alpha$ is a solution
then $q^4vv^*u^{-1}(u^*)^{-1} \alpha^{-1}$ is also a solution
and there is no further solution.
\end{note}

\begin{proposition}
\label{prop2:anyp}
Let $d$ denote a feasible integer and 
let 
$P \in  \K\lbrack \lambda 
\rbrack$ 
denote a monic polynomial 
of degree $d$. Then there exists a 
standard 
$U\sb q(\widehat{\mathfrak{sl}}\sb 2)$-module
$V$ of diameter $d$ such that 
$P_V=P$.
\end{proposition}
\noindent {\it Proof:} 
Since $\K$ is algebraically closed
there exist scalars $\lbrace r_i\rbrace_{i=1}^d$
in $\K$
such that $P= \prod_{i=1}^d (\lambda - r_i)$. By
Lemma
\ref{lem2:allp}, for $1 \leq i \leq d$ there exists
a nonzero $\alpha_i \in \K$
such that $\lambda - r_i$ is the
Drinfel'd polynomial for $V(\alpha_i)$.
Define
the $U\sb q(\widehat{\mathfrak{sl}}\sb 2)$-module
$V = \otimes_{i=1}^d V( \alpha_i)$.
By construction $V$ is standard with diameter $d$.
Also $P_V= \prod_{i=1}^d (\lambda - r_i)$
by Proposition
\ref{thm2:pvpwmain},
 so $P_V=P$. The result follows.
\hfill $\Box$ \\

\noindent {\it Proof of Proposition 
\ref{cor:splitany}}.
For $0 \leq i \leq d$ define
\begin{eqnarray}
\sigma_i = \frac{\zeta_i}{
(q-q^{-1})^2 
(q^2-q^{-2})^2 
\cdots 
(q^i-q^{-i})^2}.
\label{eq:sigzeta}
\end{eqnarray}
Define a polynomial
 $P \in  \K\lbrack \lambda \rbrack$
by 
\begin{eqnarray}
P =
(-1)^d \sum_{i=0}^d 
 \sigma_{d-i} f_0 f_1 \cdots f_{i-1}
\label{eq:pfix}
\end{eqnarray}
where the $f_j$ are from
Definition \ref{def:pi}.
Observe that $P$ is monic with degree $d$.
By 
Proposition
\ref{prop2:anyp}
there exists a standard
$U\sb q(\widehat{\mathfrak{sl}}\sb 2)$-module
of diameter $d$ such that
$P_V=P$.
Comparing
(\ref{eq:drinp}),
(\ref{eq:pfix})
 we find
the sequence
$\lbrace \sigma_i \rbrace_{i=0}^d$ from
(\ref{eq:sigzeta})
is the normalized split sequence
for $V$.
Comparing 
(\ref{eq:splitnorm}),
(\ref{eq:sigzeta})
we find
$\lbrace \zeta_i \rbrace_{i=0}^d$ is the split
sequence for $V$.
\hfill $\Box$ \\

\noindent We finish this section with some 
formulae for later use.

\medskip
\noindent Let $V$ denote a standard 
$U\sb q(\widehat{\mathfrak{sl}}\sb 2)$-module 
with feasible diameter $d \geq 1$, and let
$\lbrace U_i\rbrace_{i=0}^d$ denote the corresponding
weight space decomposition.
Let $\zeta_1$ denote term one of
the split sequence for $V$
and recall from Definition
\ref{def:splituq} that $\zeta_1$ is the
eigenvalue of $LR$ on $U_0$.
By Lemma
\ref{lem:rlmove} the space $U_d$ is invariant
under $RL$; let
$\zeta^\times_1$ denote the corresponding eigenvalue.  

\begin{lemma}
\label{lem:zetaone}
Let $V=\otimes_{i=1}^d V(\alpha_i)$ denote
a standard $U\sb q(\widehat{\mathfrak{sl}}\sb 2)$-module 
with feasible diameter $d \geq 1$. 
Then 
the following 
{\rm (i)}--{\rm (iii)} hold.
\begin{itemize}
\item[\rm (i)]
$\zeta_1
= uu^*q^{-1}\sum_{i=1}^d \alpha_i
+ 
vv^*q^3\sum_{i=1}^d \alpha^{-1}_i
-
(q-q^{-1})(q^d-q^{-d})(bb^*q^{1-d}+cc^*q^{d-1})$;
\item[\rm (ii)]
$\zeta^\times_1 = uu^*q^{-1}\sum_{i=1}^d \alpha_i
+ 
vv^*q^3\sum_{i=1}^d \alpha^{-1}_i
-
(q-q^{-1})(q^d-q^{-d})(bb^*q^{d-1}+cc^*q^{1-d})$;
\item[\rm (iii)]
$\zeta_1-\zeta^\times_1 =
(q-q^{-1})(q^{d-1}-q^{1-d})(q^d-q^{-d})(bb^*-cc^*).$
\end{itemize}
\end{lemma} 
\noindent {\it Proof:} 
Parts (i), (ii) are routinely obtained using
(\ref{eq:kus}), 
(\ref{eq:e1maction})--(\ref{eq:e0paction})
and Definition
\ref{def:RL}.
Part (iii) follows from (i), (ii).
\hfill $\Box$ \\

\section{The elements $A,A^*$ of 
$U\sb q(\widehat{\mathfrak{sl}}\sb 2)$}

\noindent From now until the end of Section 9
we adopt the following assumption.

\begin{assumption}
\label{not:settup}
\label{not:RL}
\rm
Assume the field $\K$ is algebraically closed.
We fix a $q$-Racah sequence
$(\lbrace \theta_i\rbrace_{i=0}^d;
\lbrace \theta^*_i\rbrace_{i=0}^d)$
 of scalars in $\K$.
We fix 
$q, a,b,c,a^*,b^*,c^*$ 
that satisfy (\ref{eq:cond1})--(\ref{eq:cond4})
and fix $R,L \in  
U\sb q(\widehat{\mathfrak{sl}}\sb 2)$
as in Definition
\ref{def:RL}.
\end{assumption}

\noindent In this section we define the elements
$A,A^* \in 
U\sb q(\widehat{\mathfrak{sl}}\sb 2)$ and investigate their 
properties.

\begin{definition}
\label{def:aas}
\rm
With reference to 
Assumption \ref{not:settup}
we define
\begin{eqnarray}
A &=& a 1 + bK_0 + cK_1 + R,
\label{eq:adef}
\\
A^* &=& a^* 1 + b^*K_0 + c^*K_1 + L.
\label{eq:asdef}
\end{eqnarray}
\end{definition}

\noindent We are going to show that $A,A^*$ satisfy
a pair of tridiagonal relations.
We now introduce the parameters involved in those relations.

\begin{definition}
\label{def:gamdel}
\rm
Define
$\beta=q^2+q^{-2}$ and 
\begin{eqnarray}
\label{eq:gamdel}
&&\gamma \;=\;  -a(q-q^{-1})^2,
\qquad \qquad 
\varrho \;=\; a^2(q-q^{-1})^2-bc(q^2-q^{-2})^2,
\\
\label{eq:gamdels}
&&\gamma^* \;=\; 
 -a^*(q-q^{-1})^2,
\quad \qquad 
\varrho^* = a^{*2}(q-q^{-1})^2-b^*c^*(q^2-q^{-2})^2.
\end{eqnarray}
\end{definition}

\noindent We mention one significance of the
parameters in Definition
\ref{def:gamdel}.

\begin{lemma}
\label{lem:gamdelmeaning}
The
 following {\rm (i)}--{\rm (iv)} hold.
\begin{itemize}
\item[\rm (i)]
$\gamma=\theta_{i-1}-\beta \theta_i + \theta_{i+1} \qquad (1 \leq i \leq d-1)$;
\item[\rm (ii)]
$\gamma^*=\theta^*_{i-1}-\beta \theta^*_i + \theta^*_{i+1} 
\qquad (1 \leq i \leq d-1)$;
\item[\rm (iii)]
$\varrho= \theta^2_{i-1}-\beta \theta_{i-1}\theta_i+\theta^2_i
-\gamma (\theta_{i-1}+\theta_i) \qquad (1 \leq i \leq d)$;
\item[\rm (iv)]
$\varrho^*= \theta^{*2}_{i-1}-\beta \theta^*_{i-1}\theta^*_i+\theta^{*2}_i
-\gamma^* (\theta^*_{i-1}+\theta^*_i) \qquad (1 \leq i \leq d)$.
\end{itemize}
\end{lemma}
\noindent {\it Proof:} 
Routine verification using
(\ref{eq:cond1}), (\ref{eq:cond2}) and Definition
\ref{def:gamdel}.
\hfill $\Box $ \\

\begin{proposition}
\label{lem:tdrel}
In
$U_q(L(\mathfrak{sl}\sb 2))$,
\begin{eqnarray*}
A^3A^*- \lbrack 3\rbrack_q A^2A^*A+
\lbrack 3\rbrack_q AA^*A^2-
A^*A^3 &=& 
\gamma (A^2A^*-A^*A^2)+\varrho (AA^*-A^*A),
\\
A^{*3}A- \lbrack 3\rbrack_q A^{*2}AA^*+
\lbrack 3\rbrack_q A^*AA^{*2}-
AA^{*3} &=& 
\gamma^* (A^{*2}A-AA^{*2})+\varrho^* (A^*A-AA^*),
\end{eqnarray*}
where $\gamma, \gamma^*, \varrho, \varrho^*$ are from
Definition
\ref{def:gamdel}.
\end{proposition}
\noindent {\it Proof:} 
Routine verification using
Lemma \ref{lem:RLK},
Lemma
\ref{lem:rlk},
and
Definition \ref{def:aas}.
\hfill $\Box $ \\

\begin{note}
\rm The above equations are called the {\it tridiagonal relations}
  \cite{qSerre}.
\end{note}

\section{The action of $A,A^*$ on a standard 
$U\sb q(\widehat{\mathfrak{sl}}\sb 2)$-module}

\noindent Throughout this section Assumption
\ref{not:settup}
remains in effect. We fix a standard 
$U\sb q(\widehat{\mathfrak{sl}}\sb 2)$-module
$V$ with diameter $d$ and let $\lbrace U_i\rbrace_{i=0}^d
$ denote the corresponding weight space decomposition.
Let $\lbrace\zeta_i\rbrace_{i=0}^d$ denote the
split sequence of $V$.

\medskip
\noindent In this section we describe the action of $A,A^*$ on $V$.

\begin{lemma}
\label{lem:qavoid}
We have
$q^{2i}\not=1$
for $1 \leq i \leq d$. In other words
$d$ is feasible with respect to $q$.
\end{lemma}
\noindent {\it Proof:} 
Immediate from Lemma
\ref{lem:clar}.
\hfill $\Box$ \\

\begin{lemma}
\label{lem:acton}
The following 
 hold
for  $0 \leq i \leq d$:
\begin{itemize}
\item[\rm (i)]
The element $R$ acts on $U_i$ as $A-\theta_i 1$.
\item[\rm (ii)]
The element $L$ acts on $U_i$ as $A^*-\theta^*_i 1$.
\end{itemize}
\end{lemma}
\noindent {\it Proof:} 
Immediate from
Definition
\ref{def:aas}
and 
(\ref{eq:nnkmove}).
\hfill $\Box$ \\

\begin{lemma} 
\label{lem:acton2}
The following hold 
for $0 \leq i \leq d$:
\begin{enumerate}
\item[{\rm (i)}]  
$(A-\theta_i 1)U_i \subseteq U_{i+1}$,
\item[{\rm (ii)}]  
$(A^*-\theta^*_i 1)U_i \subseteq U_{i-1}$.
\end{enumerate}
\end{lemma}
\noindent {\it Proof:} 
(i) 
Combine
Lemma
\ref{lem:acton}(i) with the inclusion on
the left in
(\ref{eq:rlmove}).
\\
\noindent (ii)
Combine Lemma
\ref{lem:acton}(ii) with the inclusion on
the right in
(\ref{eq:rlmove}).
\hfill $\Box $ \\

\begin{lemma}
\label{lem:aeig}
The element $A$ (resp. $A^*$) is diagonalizable
on $V$ with eigenvalues 
$\lbrace \theta_i \rbrace_{i=0}^d$
(resp.
$\lbrace \theta^*_i \rbrace_{i=0}^d$).
Moreover for $0 \leq i \leq d$
the dimension of the eigenspace for $A$ (resp. $A^*$)
associated with $\theta_i$ (resp. $\theta^*_i$)
is equal to 
$\binom{d}{i}$.
\end{lemma}
\noindent {\it Proof:} 
We first display the eigenvalues of $A$.
Recall that
$\lbrace \theta_i \rbrace_{i=0}^d$ are mutually
distinct, and 
$U_i$ has dimension $\binom{d}{i}$ for
$0 \leq i \leq d$.
By Lemma \ref{lem:acton2}(i)
we see that, 
with respect to an appropriate basis for
$V$, $A$ is represented by a lower triangular
matrix that has diagonal entries
$\lbrace \theta_i\rbrace_{i=0}^d$, with $\theta_i$ appearing
$\binom{d}{i}$  times for $0 \leq i \leq d$.
Hence for $0 \leq i \leq d$ the scalar
$\theta_i$ is a root of the
characteristic polynomial of $A$ with multiplicity
$\binom{d}{i}$.
We now show $A$ is diagonalizable.
To do this we show that the minimal polynomial of
$A$ has distinct roots.
By Lemma \ref{lem:acton2}(i)
we find $\prod_{i=0}^d(A-\theta_i 1)$ vanishes on $V$.
By this and since
$\lbrace \theta_i \rbrace_{i=0}^d$ are distinct we see that
the minimal polynomial of $A$ has distinct roots.
Therefore $A$ is diagonalizable.
We have now proved our assertions concerning $A$; 
our assertions concerning $A^*$ are similarly proved.
\hfill $\Box $ \\

\noindent At this point it is convenient to introduce
the primitive idempotents for $A$ and $A^*$.

\begin{definition}
\label{def:edef}
\rm
For $0 \leq i \leq d$ we define the following elements
in 
$U\sb q(\widehat{\mathfrak{sl}}\sb 2)$:
\begin{eqnarray}
\label{eq:eiesi}
  E_i=\prod_{\stackrel{0 \leq j \leq d}{j \neq i}}
          \frac{A-\theta_j 1}{\theta_i-\theta_j},
\qquad
  E^*_i=\prod_{\stackrel{0 \leq j \leq d}{j \neq i}}
          \frac{A^*-\theta^*_j 1}{\theta^*_i-\theta^*_j}.
\end{eqnarray}
We observe that $E_i$ (resp. $E^*_i$) 
acts on $V$ as the primitive idempotent of $A$ (resp. $A^*$) associated with
the eigenvalue $\theta_i$ (resp. $\theta^*_i$).
In particular $E_iV$ (resp. $E^*_iV$) is the eigenspace of
$A$ (resp. $A^*$) on $V$ associated with
the eigenvalue $\theta_i$ (resp. $\theta^*_i$).
\end{definition}

\begin{lemma}
\label{lem:eidim}
Each of $E_iV$ and $E^*_iV$ has dimension 
$\binom{d}{i}$
for $0 \leq i \leq d$.
\end{lemma}
\noindent {\it Proof:} 
Immediate from
Definition \ref{def:edef}
and the last sentence in
Lemma
\ref{lem:aeig}.
\hfill $\Box $ \\

\begin{lemma}
\label{lem:dsum}
The following hold for
$0 \leq i \leq d$:
\begin{enumerate}
\item[{\rm (i)}]  
$E_iV+ \cdots + E_dV = U_i+ \cdots+ U_d$, 
\item[{\rm (ii)}]  
$E^*_0V+ \cdots + E^*_iV  = U_{0}+ \cdots+ U_i$. 
\end{enumerate}
\end{lemma}
\noindent {\it Proof:} 
\noindent (i) 
Let $X_i=\sum_{j=i}^d U_j$
and 
$X'_i=\sum_{j=i}^d E_jV$.
We show $X_i=X'_i$.
Define $T_i = \prod_{j=i}^d (A-\theta_j 1)$.
Then $X'_i = \lbrace v \in V\,|\,T_iv=0\rbrace$,
and $T_iX_i=0$ by
Lemma \ref{lem:acton2}(i),
so $X_i \subseteq X'_i$.
Now define $S_i= \prod_{j=0}^{i-1}(A-\theta_j 1)$.
Observe that $S_iV=X'_i$, and
 $S_iV\subseteq X_i$ by
Lemma \ref{lem:acton2}(i),
so $X'_i \subseteq X_i$.
By these comments $X_i=X'_i$. 
\\
\noindent (ii) Similar to the proof of (i) above.
\hfill $\Box $ \\

\begin{lemma}
\label{lem:a0} The following
{\rm (i)}, {\rm (ii)} hold on $V$ provided that $d\geq 1$.
\begin{itemize}
\item[\rm (i)]
 $E^*_0AE^*_0=a_0E^*_0$
where
\begin{eqnarray}
a_0 = \theta_0 + \zeta_1 (\theta^*_0-\theta^*_1)^{-1}.
\label{eq:a0}
\end{eqnarray}
\item[\rm (ii)]
 $E_dA^*E_d=a^*_dE_d$
where
\begin{eqnarray}
a^*_d = \theta^*_1 - \frac{\zeta_1 +(\theta^*_0-\theta^*_1)(\theta_0
-\theta_{d-1})}{\theta_{d-1}-\theta_d}.
\label{eq:ads}
\end{eqnarray}
\end{itemize}
\end{lemma}
\noindent {\it Proof:} 
(i) The expression $E^*_0-1$ is zero on $E^*_0V$
by Definition
\ref{def:edef}, and $E^*_0V=U_0$
by Lemma
\ref{lem:dsum}(ii). Therefore 
it suffices to show that 
$E^*_0(A-a_0 1)$ is zero on $U_0$. 
By 
Definition 
\ref{def:splituq}
$LR-\zeta_1 1$ is zero on $U_0$.
By  
Lemma \ref{lem:acton}
and Lemma
\ref{lem:acton2},
$LR=(A^*-\theta^*_1 1)(A-\theta_0 1)$
on $U_0$. By Definition
\ref{def:edef},
$E^*_0A^*=\theta^*_0E^*_0$ on $V$.
We may now argue that on $U_0$,
\begin{eqnarray*}
E^*_0\zeta_1 &=&
E^*_0LR
\\
&=& E^*_0(A^*-\theta^*_1 1)(A-\theta_0 1)
\\
&=& (\theta^*_0-\theta^*_1)E^*_0(A-\theta_0 1).
\end{eqnarray*}
By this and
(\ref{eq:a0}) the expression $E^*_0(A-a_0 1)$ is zero on $U_0$.
\\
\noindent (ii)
 The expression $E_d-1$ is zero on $E_dV$
by Definition
\ref{def:edef}, and $E_dV=U_d$
by Lemma
\ref{lem:dsum}(i). Therefore 
it suffices to show that 
$E_d(A^*-a^*_d 1)$ is zero on $U_d$. 
By the paragraph above  
Lemma 
\ref{lem:zetaone},
the expression $RL-\zeta^\times_1 1$ is zero on $U_d$.
Evaluating Lemma
\ref{lem:zetaone}(iii) using
 (\ref{eq:cond1}),
(\ref{eq:cond2})
we find
\begin{eqnarray}
\zeta_1 - \zeta^\times_1 =
(\theta^*_1-\theta^*_d)
(\theta_{d-1}-\theta_d)
-
(\theta^*_0-\theta^*_1)
(\theta_0-\theta_{d-1}).
\label{eq:rlud}
\end{eqnarray}
By  
Lemma \ref{lem:acton}
and Lemma
\ref{lem:acton2},
$RL=(A-\theta_{d-1} 1)(A^*-\theta^*_d 1)$
on $U_d$. By
Definition
\ref{def:edef},
$E_dA=\theta_dE_d$ on $V$.
We may now argue that on $U_d$,
\begin{eqnarray*}
E_d \zeta^\times_1 &=&
E_d RL
\\
&=&
E_d(A-\theta_{d-1} 1)(A^*-\theta^*_d 1)
\\
&=&
(\theta_d-\theta_{d-1})E_d(A^*-\theta^*_d 1).
\end{eqnarray*}
By this and
 (\ref{eq:ads}),
(\ref{eq:rlud})
 the expression $E_d(A^*-a^*_d 1)$ is zero on $U_d$.
\hfill $\Box $ \\

\begin{lemma}
\label{lem:quasitd}
For $0 \leq i,j\leq d$ the following 
{\rm (i)},  
{\rm (ii)} hold on $V$.    
\begin{enumerate} 
\item[{\rm (i)}]  $E_iA^*E_j=0$  if  $|i-j|>1$;
\item[{\rm (ii)}]  $E^*_iAE^*_j=0$ if $|i-j|>1$.
\end{enumerate}
\end{lemma}
\noindent {\it Proof:} 
Assume $d\geq 2$; otherwise there is nothing to prove.
We first show that $E_iA^*E_j=0$ on $V$ for
$1 < |i-j|<d$.
For notational convenience define
a two variable polynomial
\begin{eqnarray}
p(\lambda, \mu)= \lambda^2-\beta \lambda \mu + \mu^2-\gamma (\lambda+\mu)
-\varrho,
\label{eq:defp}
\end{eqnarray}
where $\beta$,
$\gamma$, $\varrho$ are from
Definition \ref{def:gamdel}.
In the first equation of Proposition
\ref{lem:tdrel} 
we multiply each term on the left
by $E_i$ and the right by $E_j$.
We simplify the result using
the fact that $E_iA=\theta_iE_i$
and 
$AE_j=\theta_jE_j$ on $V$.
After a brief calculation this shows
\begin{eqnarray}
\label{eq:wantthis}
0 = E_iA^*E_j(\theta_i-\theta_j)p(\theta_i,\theta_j)
\end{eqnarray}
on $V$. 
We claim that in (\ref{eq:wantthis}) the
coefficient of $E_iA^*E_j$ is nonzero.
Of course $\theta_i-\theta_j\not=0$ since
$\theta_0, \ldots, \theta_d$ are mutually distinct. 
We now show
 $p(\theta_i,\theta_j)\not=0$.
Since $|i-j|<d$ we have 
$1 \leq i \leq d-1$
or $1 \leq j \leq d-1$. We may assume $1 \leq i \leq d-1$
since $p(\theta_i,\theta_j)=
 p(\theta_j,\theta_i)$.
By Lemma
\ref{lem:gamdelmeaning}(iii) 
we have $p(\theta_i,\theta_{i-1})=0$
and 
 $p(\theta_i,\theta_{i+1})=0$. The expression
$p(\theta_i, \mu)$ is a quadratic polynomial in $\mu$ with 
roots $\theta_{i-1}$, 
$\theta_{i+1}$.
Since $|i-j|>1$ we have
$\theta_j \not=\theta_{i-1}$
and
$\theta_j \not=\theta_{i+1}$.
Therefore $p(\theta_i,\theta_j)\not=0$ as desired.
Now in (\ref{eq:wantthis}) the coefficient of 
$E_iA^*E_j$ is nonzero, so
$E_iA^*E_j=0$ on $V$.
Swapping the roles of $A,A^*$ in the above argument,
we similarly find
$E^*_iAE^*_j=0$ on $V$ for $1 < |i-j|<d$.
Next we show that $E_0A^*E_d=0$ on $V$.
Observe 
$(A^*-\theta^*_d 1)U_d \subseteq U_{d-1}$
by Lemma
\ref{lem:acton2}(ii), so
$A^*U_d \subseteq U_{d-1}+U_d$.
By Lemma
\ref{lem:dsum}(i),
 $E_dV=U_d$ and
$E_{d-1}V+E_dV=U_{d-1}+U_d$.
By these comments
$A^*E_dV\subseteq E_{d-1}V+E_dV$. We assume
$d\geq 2$ so $E_0$ vanishes 
on 
$E_{d-1}V+E_dV$ and 
therefore $E_0A^*E_d=0$ on $V$.
 Next we show that $E^*_dAE^*_0=0$ on $V$.
Observe
$(A-\theta_0 1)U_0 \subseteq U_1$ by
Lemma
\ref{lem:acton2}(i), so
$AU_0 \subseteq U_0+U_1$.
By Lemma
\ref{lem:dsum}(ii),
 $E^*_0V=U_0$ and
$E^*_0V+E^*_1V=U_0+U_1$.
By these comments
$AE^*_0V\subseteq E^*_0V+E^*_1V$. We assume
$d\geq 2$ so $E^*_d$ vanishes 
on 
$E^*_0V+E^*_1V$ and
therefore $E^*_dAE^*_0=0$ on $V$.
Next we show that $E_dA^*E_0=0$ on $V$.
Since $V=\sum_{i=0}^d U_i$
it suffices to show
that
$E_dA^*E_0=0$ on $U_i$ for $0 \leq i \leq d$.
Observe 
$E_dA^*E_0=0$ on $U_i$ for $1 \leq i \leq d$, since
$\sum_{i=1}^d U_i=\sum_{i=1}^d E_iV$ by
Lemma
\ref{lem:dsum}(i) and since $E_0$ vanishes on
$E_iV$ for $1 \leq i \leq d$.
To show $E_dA^*E_0=0$ on $U_0$, 
recall $U_0=E^*_0V$ 
so it suffices to show
that $E_dA^*E_0E^*_0=0$ on $V$.
By Definition \ref{def:edef},
$1=\sum_{i=0}^d E^*_i$ on $V$. In this equation
we multiply
each term on the right by 
$AE^*_0$. We evaluate the result using Lemma
\ref{lem:a0}(i) and the fact that
$E^*_iAE^*_0=0$ on $V$ for $2 \leq i \leq d$.
This yields
\begin{eqnarray}
\label{eq:p3}
AE^*_0 = a_0 E^*_0 + E^*_1AE^*_0
\end{eqnarray}
on $V$.
In (\ref{eq:p3}) we multiply each term on the left by $A^*$
to find
\begin{eqnarray}
A^*AE^*_0 = a_0\theta^*_0E^*_0+ \theta^*_1E^*_1AE^*_0
\label{eq:p4}
\end{eqnarray}
on $V$.
By Definition
\ref{def:edef},
$1=\sum_{i=0}^d E_i$ on $V$.
In this equation we multiply each term on the left
by $E_dA^*$. We evaluate the result using
Lemma
\ref{lem:a0}(ii) and the fact that
$E_dA^*E_i=0$ on $V$ for $1 \leq i \leq d-2$; this yields
\begin{eqnarray}
E_dA^*=
E_dA^*E_0+
E_dA^*E_{d-1}
+ a^*_dE_d
\label{eq:p1}
\end{eqnarray}
on $V$.
In 
(\ref{eq:p1}) we multiply each term
on the right by $A$ to find
\begin{eqnarray}
E_dA^*A=
\theta_0 E_dA^*E_0+
\theta_{d-1}E_dA^*E_{d-1}
+ \theta_d a^*_dE_d
\label{eq:p2}
\end{eqnarray}
on $V$.
 Consider the equation which is
 $\theta^*_1E_d$ times
(\ref{eq:p3})  minus
 $E_d$ times 
(\ref{eq:p4})
minus
(\ref{eq:p1}) times $\theta_{d-1}E^*_0$
plus 
(\ref{eq:p2}) times $E^*_0$.
We simplify this equation
using the fact that $A^*E^*_0=\theta^*_0E^*_0$
and 
$E_dA=\theta_dE_d$ on $V$.
The calculation shows
that on $V$ the expression 
$(\theta_0-\theta_{d-1})E_dA^*E_0E^*_0$ coincides with
 $E_dE^*_0$ times
\begin{eqnarray}
\label{eq:fincoeff}
(\theta^*_0-\theta^*_1)a_0
+
(\theta_{d-1}-\theta_d)a^*_d
+ \theta_d\theta^*_1-\theta_{d-1}\theta^*_0.
\end{eqnarray}
Note that $\theta_0-\theta_{d-1}$ is nonzero since
$d\geq 2$.
By
(\ref{eq:a0}), 
(\ref{eq:ads})  
the expression (\ref{eq:fincoeff}) is zero.
Therefore 
$E_dA^*E_0E^*_0=0$ on $V$ and
hence
$E_dA^*E_0=0$ on $V$.
Next we show that
$E^*_0AE^*_d=0$ on $V$.
Since $V=\sum_{i=0}^d U_i$ 
it suffices to show
that
$E^*_0AE^*_d=0$ on $U_i$ for $0 \leq i \leq d$.
Observe 
$E^*_0AE^*_d=0$ on $U_i$ for $0 \leq i \leq d-1$, since
$\sum_{i=0}^{d-1} U_i=\sum_{i=0}^{d-1} E^*_iV$ by
Lemma
\ref{lem:dsum}(ii) and since $E^*_d$ vanishes on
$E^*_iV$ for $0 \leq i \leq d-1$.
To show $E^*_0AE^*_d=0$ on $U_d$, 
recall $U_d=E_dV$ 
so it suffices to show
that
$E^*_0AE^*_dE_d=0$ on $V$.
We mentioned earlier 
that $1=\sum_{i=0}^d E_i$ on $V$. In this equation
we multiply
each term on the right by 
$A^*E_d$. We evaluate the result using Lemma
\ref{lem:a0}(ii) and the fact that
$E_iA^*E_d=0$ on $V$ for $0 \leq i \leq d-2$.
This yields
\begin{eqnarray}
\label{eq:p3s}
A^*E_d =
E_{d-1}A^*E_d
+
 a^*_d E_d 
\end{eqnarray}
on $V$.
In 
(\ref{eq:p3s}) we multiply each term
on the left by $A$ to find
\begin{eqnarray}
\label{eq:p4s}
AA^*E_d =
\theta_{d-1}E_{d-1}A^*E_d
+
 \theta_d a^*_d E_d 
\end{eqnarray}
on $V$.
We mentioned earlier that
$1=\sum_{i=0}^d E^*_i$ on $V$.
In this equation we multiply each term on the left
by $E^*_0A$. We evaluate the result using
Lemma
\ref{lem:a0}(i) and the fact that
$E^*_0AE^*_i=0$ on $V$ for $2 \leq i \leq d-1$; this yields
\begin{eqnarray}
E^*_0A=
a_0E^*_0
+
E^*_0AE^*_1+
E^*_0AE^*_d
\label{eq:p1s}
\end{eqnarray}
on $V$.
In 
(\ref{eq:p1s}) we multiply each term on the right by $A^*$ to find
\begin{eqnarray}
E^*_0AA^*=
\theta^*_0a_0E^*_0
+
\theta^*_1E^*_0AE^*_1+
\theta^*_d E^*_0AE^*_d
\label{eq:p2s}
\end{eqnarray}
on $V$.
Consider the equation which is
 $\theta_{d-1}E^*_0$ times
(\ref{eq:p3s})
  minus
 $E^*_0$ times 
(\ref{eq:p4s})
minus
(\ref{eq:p1s}) times $\theta^*_1 E_d$  
plus 
(\ref{eq:p2s}) times $E_d$.
We simplify this equation
using the fact that $AE_d=\theta_dE_d$
and 
$E^*_0A^*=\theta^*_0E^*_0$ on $V$.
The calculation shows
that on $V$ the expression 
$(\theta^*_1-\theta^*_d)E^*_0AE^*_dE_d$ coincides with
 $E^*_0E_d$ times
(\ref{eq:fincoeff}).
The scalar 
$\theta^*_1-\theta^*_d$ is nonzero
and we already showed that
(\ref{eq:fincoeff}) is zero, so
$E^*_0AE^*_dE_d=0$ on $V$ and
hence
$E^*_0AE^*_d=0$ on $V$.
\hfill $\Box $ \\

\begin{lemma}     \label{lem:tauc}
For $0 \leq i \leq d$ the following holds on $V$:
\begin{equation}     \label{eq:tauc}
E^*_0\tau_i(A)E^*_0 =
 \frac{\zeta_iE^*_0}
      {(\theta^*_0-\theta^*_1)(\theta^*_0-\theta^*_2)\cdots(\theta^*_0-\theta^*_i)}.
\end{equation}
\end{lemma}

\noindent {\it Proof:} 
It suffices to show that
\begin{eqnarray}
(\theta^*_0-\theta^*_1)
(\theta^*_0-\theta^*_2) \cdots
(\theta^*_0-\theta^*_i)E^*_0\tau_i(A)-\zeta_i 1
\label{eq:wantzero}
\end{eqnarray}
is zero on $E^*_0V$. We pick $w \in E^*_0V$
and show 
(\ref{eq:wantzero}) is zero at $w$.
Setting $i=0$ in 
Lemma
\ref{lem:dsum}(ii) we find
$E^*_0V=U_0$.
By 
Definition
\ref{def:splituq}
 we find
$L^iR^i-\zeta_i 1$ is zero on $U_0$. By these comments
$L^iR^i w=\zeta_i w$.
Using 
Lemma \ref{lem:acton}
and \ref{lem:acton2},
\begin{eqnarray*}
L^iR^i=
(A^*-\theta^*_1 1)
(A^*-\theta^*_2 1) \cdots
(A^*-\theta^*_i 1)\tau_i(A)
\end{eqnarray*}
on $U_0$. Therefore
\begin{eqnarray*}
(A^*-\theta^*_1 1)
(A^*-\theta^*_2 1) \cdots
(A^*-\theta^*_i 1)\tau_i(A)w=\zeta_i w.
\end{eqnarray*}
In this equation we apply $E^*_0$ to both sides
and use
$E^*_0A^*=\theta^*_0E^*_0$,
$E^*_0w=w$ to find 
(\ref{eq:wantzero}) is zero at $w$. The result follows.
\hfill $\Box$ \\

\section{The proof of
Theorem
 \ref{thm:mainth}}

\noindent Throughout this section 
we adopt the following assumption.

\begin{assumption}
\label{ass:final}
\rm Assume the field $\K$ is algebraically closed.
We fix a $q$-Racah sequence
$(\lbrace \theta_i\rbrace_{i=0}^d;
\lbrace \theta^*_i\rbrace_{i=0}^d)$
 of scalars in $\K$, and
a sequence $\lbrace \zeta_i\rbrace_{i=0}^d$
of scalars in $\K$ that satisfy
condition (ii) of Theorem
 \ref{thm:mainth}.
\end{assumption}

\noindent With reference to  Assumption \ref{ass:final},
our goal in this section is to display a TD system
over $\K$ that has parameter array
$(\lbrace \theta_i\rbrace_{i=0}^d;
\lbrace \theta^*_i\rbrace_{i=0}^d;
\lbrace \zeta_i\rbrace_{i=0}^d)$.
To this end we fix 
$q,a,b,c,a^*,b^*,c^*$ that satisfy
 (\ref{eq:cond1})--(\ref{eq:cond4}).
Using this data
we define
$R,L \in 
U\sb q(\widehat{\mathfrak{sl}}\sb 2)$
as in Definition
\ref{def:RL}, and then
$A,A^* \in 
U\sb q(\widehat{\mathfrak{sl}}\sb 2)$
as in Definition
\ref{def:aas}.
Let $\lbrace E_i\rbrace_{i=0}^d$,
$\lbrace E^*_i\rbrace_{i=0}^d$ 
be as in Definition
\ref{def:edef}.
In view of 
Proposition
\ref{cor:splitany}
and
Lemma
\ref{lem:qavoid}
we fix a standard 
$U\sb q(\widehat{\mathfrak{sl}}\sb 2)$-module
$V$ with diameter $d$ that has split sequence 
$\lbrace \zeta_i\rbrace_{i=0}^d$.

\begin{lemma}    \label{lem:triple} 
The elements $E^*_0E_0E^*_0$, $E^*_0E_dE^*_0$ are nonzero on $V$.
\end{lemma}
\noindent {\it Proof:} 
By construction $E^*_0V\not=0$.
Concerning $E^*_0E_dE^*_0$, by the equation on the left in 
(\ref{eq:eiesi})
we have  $E_d=\tau_d(A)\tau_d(\theta_d)^{-1}$.
By Lemma \ref{lem:tauc} (with $i=d$)
 $E^*_0\tau_d(A)E^*_0=\eta^*_d(\theta^*_0)^{-1}\zeta_d E^*_0$ on $V$.
Therefore 
$E^*_0E_dE^*_0 =\tau_d(\theta_d)^{-1} \eta^*_d(\theta^*_0)^{-1}\zeta_d E^*_0$
on $V$. By this and since  $\zeta_d\not=0$ we find
 $E^*_0E_dE^*_0$ is
nonzero
on $V$.
Concerning $E^*_0E_0E^*_0$, by the equation on the left in
 \eqref{eq:eiesi} 
we have  $E_0=\eta_d(A)\eta_d(\theta_0)^{-1}$.
By \cite[Proposition 5.5]{NT:mu},
$\eta_d = \sum_{i=0}^d \eta_{d-i}(\theta_0)\tau_i $.
By these comments and Lemma \ref{lem:tauc},
\[
 E^*_0E_0E^*_0 =
  E^*_0 \eta_d(\theta_0)^{-1}\eta^*_d(\theta^*_0)^{-1}
   \sum_{i=0}^d  \eta_{d-i}(\theta_0)\eta^*_{d-i}(\theta^*_0) \zeta_i
\]
on $V$.
 In the above line the sum is nonzero by \eqref{eq:ineq} so
$E^*_0E_0E^*_0$ is
 nonzero on $V$.
\hfill $\Box$ \\

\begin{definition}
\label{def:talgebra}
\rm
Let $T$ denote the subalgebra of 
$U\sb q(\widehat{\mathfrak{sl}}\sb 2)$
generated by $A,A^*$. We observe
that $T$ contains $E_i, E^*_i$ for
$0 \leq i \leq d$.
\end{definition}

\noindent Observe that $TE^*_0V$ is the $T$-submodule
of $V$ generated by $E^*_0V$. We now examine this module.

\begin{lemma}    \label{lem:110}  
Let $W$ denote a proper $T$-submodule of $TE^*_0V$.
 Then $E^*_0W=0$.
\end{lemma}
\noindent {\it Proof:} 
Suppose $E^*_0W \not=0$. 
The space $E^*_0V$ contains
$E^*_0W$ and has dimension 1, so
$E^*_0V=E^*_0W$. The space $W$ is $T$-invariant
and $E^*_0\in T$ so $E^*_0W \subseteq W$. Therefore
$E^*_0V \subseteq W$, which yields 
$TE^*_0V \subseteq W$. This contradicts
the fact that $W$ is properly contained in $TE^*_0V$.
Therefore 
 $E^*_0W =0$. 
\hfill $\Box$ \\

\begin{lemma}     \label{lem:ww}
Let $W$ and $W'$ denote proper $T$-submodules of $TE^*_0V$.
 Then $W+W'$ is a proper
$T$-submodule of $TE^*_0V$.
\end{lemma}
\noindent {\it Proof:} 
We show $W+W'\not=TE^*_0V$. 
The kernel of $E^*_0$ on $TE^*_0V$
is properly contained in $TE^*_0V$,
since $0 \not=E^*_0V \subseteq TE^*_0V$.
This kernel contains each of 
 $W,W'$ by
 Lemma \ref{lem:110}, so this kernel contains
$W+W'$.
Therefore  $W+W' \not=TE^*_0V$ 
and the result follows.
\hfill $\Box$ \\

\begin{definition}
\rm
Let $W$ denote a proper $T$-submodule of $TE^*_0V$.
 Then $W$ is called {\it maximal}
whenever $W$ is not contained in any proper $T$-submodule of $TE^*_0V$, 
besides itself.
\end{definition}

\begin{lemma}     \label{lem:120}
There exists a unique maximal proper $T$-submodule in $TE^*_0V$.
\end{lemma}
\noindent {\it Proof:} 
Concerning existence, consider  
\begin{equation}       \label{eq:sumw}
\sum_{W} W,
\end{equation}
where the sum is over all proper $T$-submodules $W$ of $TE^*_0V$.
The space \eqref{eq:sumw} is a proper $T$-submodule of $TE^*_0V$ by
Lemma \ref{lem:ww}, and since $TE^*_0V$ has finite dimension.
The $T$-submodule \eqref{eq:sumw} is maximal by the construction.
Concerning uniqueness, suppose $W$ and $W'$ are 
maximal proper $T$-submodules 
of $TE^*_0V$. By Lemma \ref{lem:ww} $W+W'$ is a proper
 $T$-submodule of $TE^*_0V$.
The space $W+W'$ contains each of $W$, $W'$ so $W+W'$ is equal to each of
$W$, $W'$ by the maximality of $W$ and $W'$. Therefore $W=W'$ and the result 
follows.
\hfill $\Box$ \\

\begin{definition}   
\rm
Let $M$ denote the maximal proper $T$-submodule of $TE^*_0V$.
Let $L$ denote the quotient $T$-module $TE^*_0V\slash M$.
By construction the $T$-module
 $L$ is nonzero, finite-dimensional  and irreducible.
\end{definition}

\begin{proposition}           \label{thm:exist}  
The sequence $(A; \{E_i\}_{i=0}^d;A^*; \{E^*_i\}_{i=0}^d)$ acts on $L$ as 
a  TD system with parameter array
 $(\{\theta_i\}_{i=0}^d; \{\theta^*_i\}_{i=0}^d; \{\zeta_i\}_{i=0}^d)$.
\end{proposition}

\noindent {\it Proof:} 
We first show that $(A; \{E_i\}_{i=0}^d; A^*; \{E^*_i\}_{i=0}^d)$ acts on
$L$ as a TD system. 
We start with a few statements that follow from the construction.
The space $L$ is a direct sum of the nonzero spaces among
$\lbrace E_iL\rbrace_{i=0}^d$ and a
 direct sum of the nonzero spaces among
$\lbrace E^*_iL\rbrace_{i=0}^d$.
For $0 \leq i \leq d$,
 $(A-\theta_i 1)E_iL=0$ 
and $(A^*-\theta^*_i 1)E^*_iL=0$.
 Using Lemma
\ref{lem:quasitd},
\begin{eqnarray}
A^*E_iL \subseteq E_{i-1}L+E_iL + E_{i+1}L
\qquad \qquad (0 \leq i \leq d),
\label{eq:l3}
\end{eqnarray}
where $E_{-1}=0$ and $E_{d+1}=0$. Moreover
\begin{eqnarray}
AE^*_iL \subseteq E^*_{i-1}L+E^*_iL + E^*_{i+1}L
\qquad \qquad (0 \leq i \leq d),
\label{eq:l3s}
\end{eqnarray}
where $E^*_{-1}=0$ and $E^*_{d+1}=0$.
Observe that $E^*_0L\not=0$ since
$M$ does not contain $E^*_0V$.
We now show 
 $E_0L\not=0$.
 Suppose $E_0L=0$.
Then $E_0TE^*_0V \subseteq M$
so
 $E_0E^*_0V \subseteq M$.
In this containment we apply $E^*_0$ to both sides and use $E^*_0M=0$ to
get $E^*_0E_0E^*_0 V=0$.
This contradicts Lemma \ref{lem:triple}
 so $E_0L\not=0$.
Next we show 
 $E_dL\not=0$.
Suppose $E_d L=0$.
Then $E_dTE^*_0V \subseteq M$
so
 $E_dE^*_0V \subseteq M$.
In this containment we apply $E^*_0$ to both sides and use $E^*_0M=0$ to
get $E^*_0E_dE^*_0 V=0$.
This contradicts Lemma \ref{lem:triple}
 so $E_dL\not=0$.
We now show $E_iL\not=0$ for $1 \leq i \leq d-1$.
Let $i$ be given and suppose 
 $E_iL=0$.
Then $E_0L+\cdots + E_{i-1}L$ is a nonzero proper 
$T$-submodule of $L$ in view of
(\ref{eq:l3}).
This contradicts the irreducibility of the $T$-module $L$.
Therefore $E_iL\not=0$ for $1 \leq i \leq d-1$.
There exists an integer $\delta$ $(0 \leq \delta \leq d)$
such that $E^*_iL\not=0$ for $0 \leq i \leq \delta$
and $E^*_{\delta+1}L=0$.
By the above comments the sequence
$(A; \{E_i\}_{i=0}^d; A^*; \{E^*_i\}_{i=0}^\delta)$ acts on
$L$ as a TD system. 
Now $d=\delta$ by the third sentence below
Note \ref{lem:convention}. 
We have shown
$(A; \{E_i\}_{i=0}^d; A^*; \{E^*_i\}_{i=0}^d)$ acts on
$L$ as a TD system which we denote by $\Phi$.
By construction $\Phi$ has eigenvalue sequence $\{\theta_i\}_{i=0}^d$
and dual eigenvalue sequence $\{\theta^*_i\}_{i=0}^d$.
By Lemma \ref{lem:tauc} and since the canonical map
 $TE^*_0V  \to L$ is a $T$-module
homomorphism, we have
\[
 E^*_0\tau_i(A)E^*_0 = 
  \frac{\zeta_i E^*_0}
       {(\theta^*_0-\theta^*_1)(\theta^*_0-\theta^*_2)
\cdots(\theta^*_0-\theta^*_i)}
     \qquad \qquad (0 \leq i \leq d)
\]
on $L$.
 By this and Definition \ref{def:split} the sequence 
$\{\zeta_i\}_{i=0}^d$ is the split sequence for $\Phi$.
By these comments $\Phi$ has parameter array
 $(\{\theta_i\}_{i=0}^d; \{\theta^*_i\}_{i=0}^d; \{\zeta_i\}_{i=0}^d)$
and the result follows.
\hfill $\Box$ \\

\noindent It is now a simple matter to prove
 Theorem
 \ref{thm:mainth}.

\medskip
\noindent {\it Proof of Theorem
 \ref{thm:mainth}}.
The implication (i)$\Rightarrow$(ii)
is proved in 
\cite[Corollary~8.3]{nomsharp}.
The implication 
(ii)$\Rightarrow$(i)
 follows
from Proposition
\ref{thm:exist}.
Now assume 
(i), (ii) hold.
Then
the last assertion of the theorem follows from
Proposition \ref{thm:isopa}.
\hfill $\Box $ \\

\section{Remarks}

\noindent In 
 this section we prove
the shape conjecture for 
the TD pairs over an algebraically closed field that have
$q$-Racah type.

\begin{proposition}
Assume the field $\K$ is algebraically closed,
and let $\lbrace \rho_i\rbrace_{i=0}^d$ denote the
shape of a TD pair over $\K$ that has $q$-Racah type.
Then $\rho_i \leq \binom{d}{i}$ for $0 \leq i \leq d$. 
\end{proposition}
\noindent {\it Proof:} 
For the TD pair in question we pick a standard ordering
of their primitive idempotents to obtain a TD system.
Without loss we may identify this
TD system with the one in   
Proposition
\ref{thm:exist}. Referring to the TD system
in 
Proposition
\ref{thm:exist},
we show 
that each of $E_iL$ and $E^*_iL$ has dimension  at most
$\binom{d}{i}$. 
The space $E_iL$ is the image of 
$E_iTE^*_0V$ under the canonical
homomorphism
$TE^*_0V \to
L$. Therefore the dimension of
$E_iL$ is at most the dimension of
$E_iTE^*_0V$.
The space $E_iTE^*_0V$ is contained in
$E_iV$ so the dimension of
 $E_iTE^*_0V$ is at most the dimension of
$E_iV$. The dimension of $E_iV$ is
$\binom{d}{i}$ 
by 
Lemma
\ref{lem:eidim}. Our conclusion for
$E_iL$ follows from the above comments. Our conclusion
for $E^*_iL$ are similarly obtained.
\hfill $\Box $ \\

\bigskip


\noindent Tatsuro Ito \hfil\break
\noindent Division of Mathematical and Physical Sciences \hfil\break
\noindent Graduate School of Natural Science and Technology\hfil\break
\noindent Kanazawa University \hfil\break
\noindent Kakuma-machi,  Kanazawa 920-1192, Japan \hfil\break
\noindent email:  {\tt tatsuro@kenroku.kanazawa-u.ac.jp}

\bigskip

\noindent Paul Terwilliger \hfil\break
\noindent Department of Mathematics \hfil\break
\noindent University of Wisconsin \hfil\break
\noindent 480 Lincoln Drive \hfil\break
\noindent Madison, WI 53706-1388 USA \hfil\break
\noindent email: {\tt terwilli@math.wisc.edu }\hfil\break

\end{document}